\documentclass[10pt, twoside]{amsart}
\usepackage{enumerate}
\usepackage{lingmacros}
\usepackage{tree-dvips}
\usepackage{amsfonts}
\usepackage{amssymb}
\usepackage{natbib}
\usepackage{mathtools}
\DeclareMathOperator*{\argmin}{arg\,min}
\usepackage[colorlinks,citecolor=blue,urlcolor=blue]{hyperref}

\usepackage{tabu}
\usepackage{longtable}
\usepackage[table]{xcolor}
\usepackage{rotating}

\usepackage{subfig}
\usepackage{dsfont}
\usepackage{pdfpages}
\usepackage{natbib}
\usepackage{setspace}
\usepackage{amsmath}
\usepackage{booktabs}
\usepackage{caption}
\usepackage{lscape}
\usepackage{graphicx}
\usepackage{amsthm}
\newtheorem{theo}{Theorem}[section]

\newtheorem{alg}[theo]{Algorithm}

\theoremstyle{definition}

\usepackage{makecell}

\definecolor{tableHeader}{RGB}{120, 120, 240}
\definecolor{tableLineOne}{RGB}{245, 245, 245}
\definecolor{tableLineTwo}{RGB}{224, 224, 224}

\newcommand{\tableHeaderStyle}{
    \rowfont{\leavevmode\color{white}\bfseries}
    \rowcolor{tableHeader}
}

\begin{document}

\title[Learning Large Scale ODE Systems]{Learning Large Scale Ordinary Differential Equation Systems}
\author[F. V. Mikkelsen]{Frederik Vissing Mikkelsen}
\address{Department of Mathematical Sciences,
University of Copenhagen,
Universitetsparken 5,
2100 Copenhagen \O,
Denmark}

\email[Corresponding author]{frm@math.ku.dk}
\author[N. R. Hansen]{Niels Richard Hansen} 
\email{Niels.R.Hansen@math.ku.dk}

\subjclass[2010]{}


\begin{abstract} Learning large scale nonlinear ordinary differential equation (ODE) systems from data is known to be computationally and statistically challenging. We present a framework together with the adaptive integral matching (AIM) algorithm for learning polynomial or rational ODE systems with a sparse network structure. The framework allows for time course data sampled from multiple environments representing e.g. different interventions or perturbations of the system. The algorithm AIM combines an initial penalised integral matching step with an adapted least squares step based on solving the ODE numerically. The R package \textit{episode} implements AIM together with several other algorithms and is available from CRAN. It is shown that AIM achieves state-of-the-art network recovery for the \emph{in silico} phosphoprotein abundance data from the eighth DREAM challenge with an AUROC of 0.74, and it is demonstrated via a range of numerical examples that AIM has good statistical properties while being computationally feasible even for large systems.  \end{abstract}

\maketitle

\noindent

\section{Introduction}
\label{sec:intro}
We consider the problem of modelling time course data sampled from a dynamical system in different environments. We model data via ordinary differential equations (ODEs), with a particular emphasis on learning the network of the system's constituents. This setting arises for instance in systems biology with biochemical reactions and molecular networks \citep{Wilkinson:2006, Oates:2012, Babtie:2014, DREAM8}, where a reaction network or a gene regulatory network may either be partially known or completely unknown. Learning such ODE networks from data is highly relevant as they provide a means for predicting downstream effects of interventions in the system.

Our main contribution is a framework and the corresponding R package \textit{episode} for learning polynomial and rational ODE systems, which is directly applicable to experimental data. Extensive numerical experiments have lead us to propose the adaptive integral matching (AIM) algorithm, though the R package includes several other learning algorithms. The framework and the learning algorithm AIM are useful when there exists little to no prior knowledge of the system in question and a fully data-driven network recognition is needed. However, the framework does also allow for incorporating prior knowledge into the estimation procedure as will be illustrated.

The paper is organised as follows. Section \ref{sec:netw} motivates the ODE network estimation problem with the small EnvZ/OmpR system. Section \ref{sec:stat} defines the statistical framework and Section \ref{sec:meth} reviews two standard approaches to parameter estimation in ODE systems: \textit{the least squares method} and \textit{the inverse collocation methods}. Then the AIM algorithm that combines both approaches is presented together with the functionality of the  R package \textit{episode} developed for this paper.  In Section \ref{sec:app} we apply our proposed method to two large scale dynamical systems: the \textit{in silico} protein phosphorylation network used in the eighth DREAM challenge \citep{DREAM8}, and a full scale model of glycolysis in \textit{Saccharomyces cerevisiae} \citep{Hynne:2001}. Finally, in Section \ref{sec:MAK} we present two extensive simulation studies that compare the performance of AIM to other methods proposed in the literature.

\section{The ODE Network Estimation Problem}
\label{sec:netw}
We illustrate the problem addressed in this paper by a simple and concrete dynamical system: the EnvZ/OmpR system. It is present among a wide range of bacteria and is particularly well studied in \textit{Escherichia coli} (\cite{Bernardini1990}, \cite{Batchelor2003}, \cite{Shinar2010}). In this system the histidine kinase $\mathrm{EnvZ}$ responds to changes in the osmolarity resulting from extracellular impermeable compounds. It responds by controlling the phosphorylation of the regulator, $\mathrm{OmpR}$, which itself proceeds to regulate the transcription of certain genes, including $ompF$ and $ompC$. These two genes act as \textit{porins} responsible for regulating the cellular diffusion across the membrane. 

The EnvZ/OmpR system is heavily studied and the whole regulation process described above is the product of numerous studies, each of which were carefully designed to isolate specific mechanisms and investigate them individually. However, various networks in systems biology are only partially understood or not even discovered yet. In this paper we do not assume that the system in question was carefully dissected and studied as a sum of local mechanisms. We simply assume that the system was observed under different perturbations or interventions and solve the network estimation problem globally. 

\begin{figure}
\begin{tabular}{cc}
\multicolumn{2}{c}{\includegraphics[width=110mm]{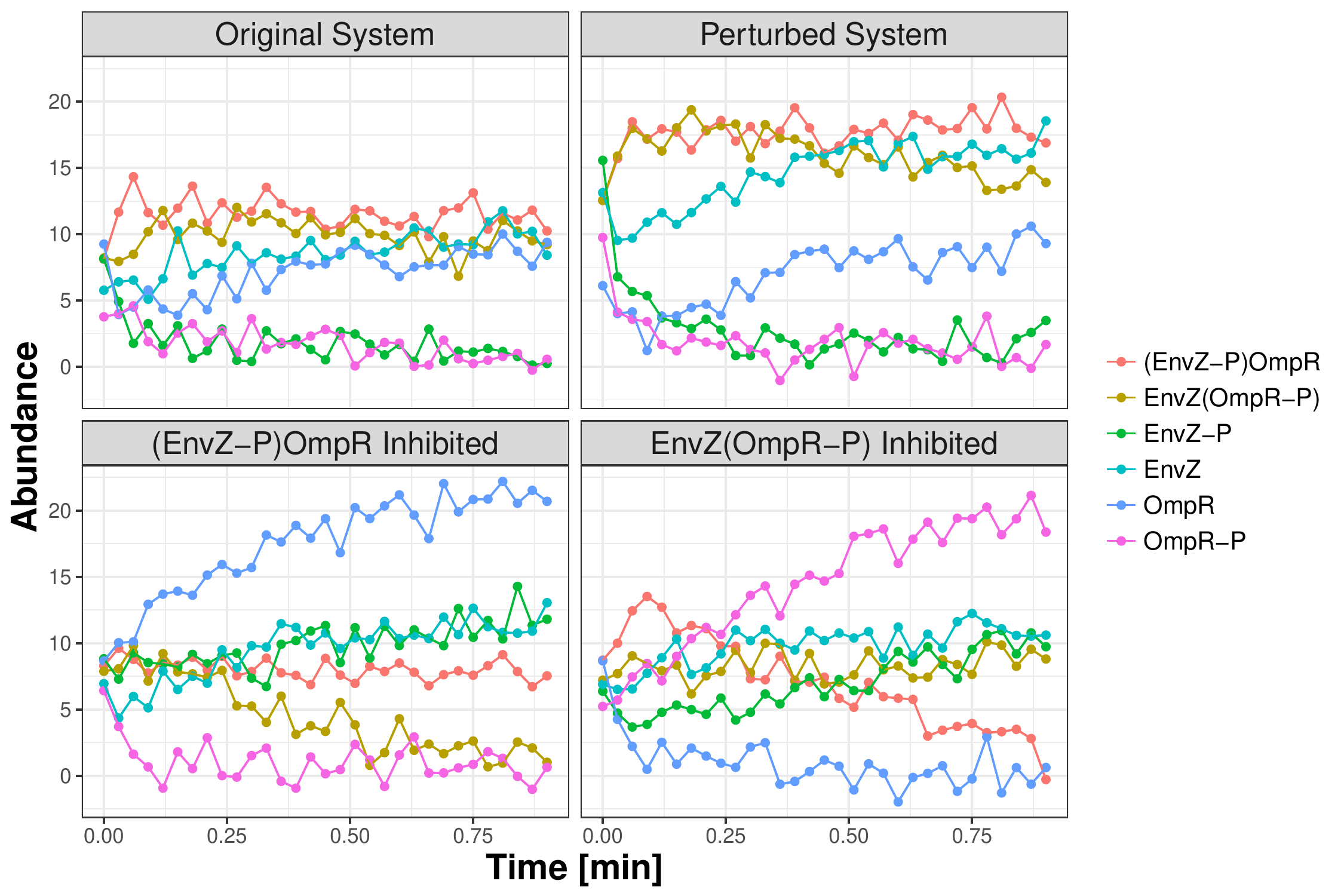}\llap{\parbox[b]{0.6in}{(a)\\\rule{0ex}{2.5in}}} }\\ \\
  \includegraphics[width=60mm]{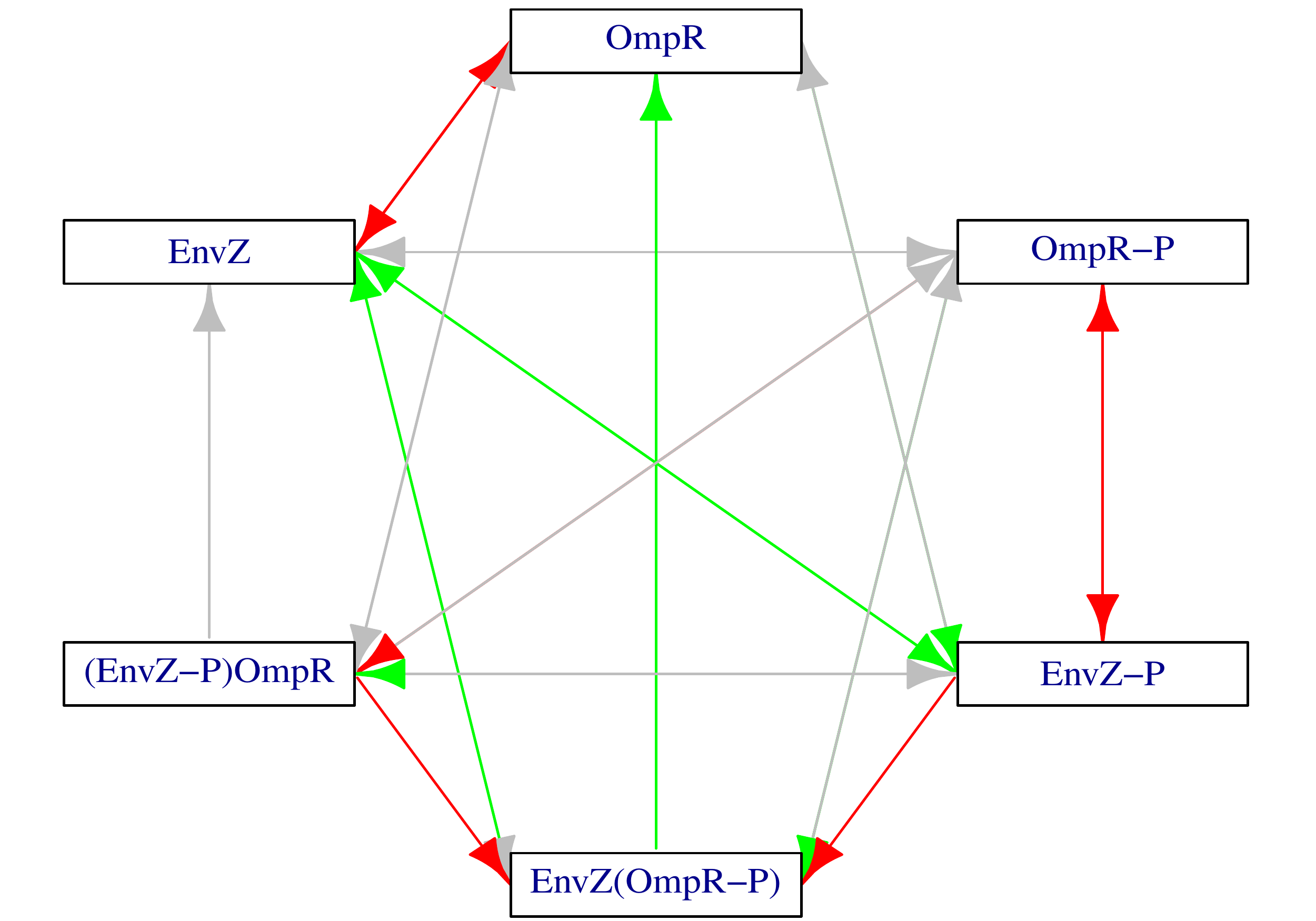}\llap{\parbox[b]{2.3in}{(b)\\\rule{0ex}{1.3in}}} & \includegraphics[width=60mm]{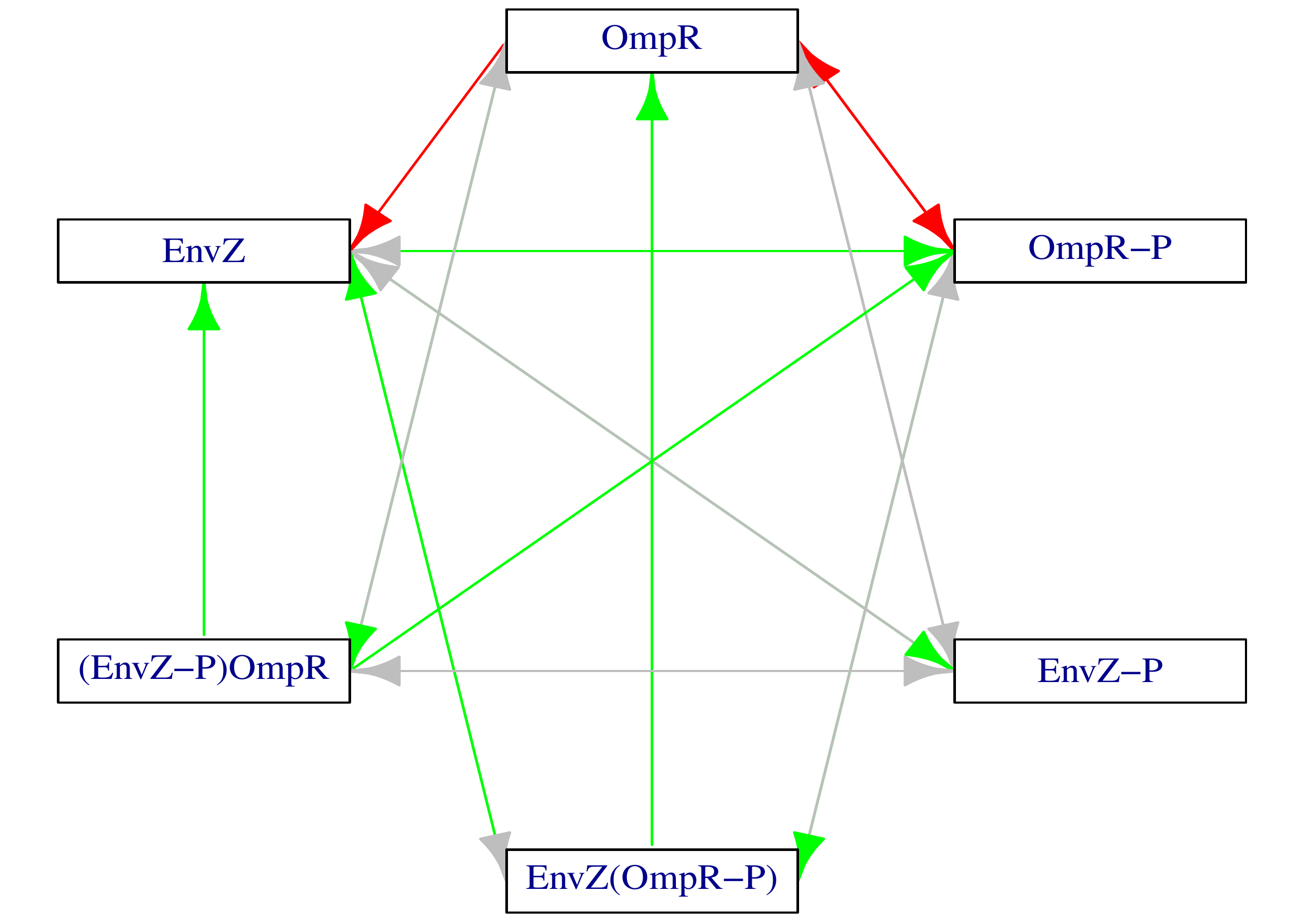}\llap{\parbox[b]{0.3in}{(c)\\\rule{0ex}{1.3in}}} \\ \\
\multicolumn{2}{c}{\includegraphics[width=60mm]{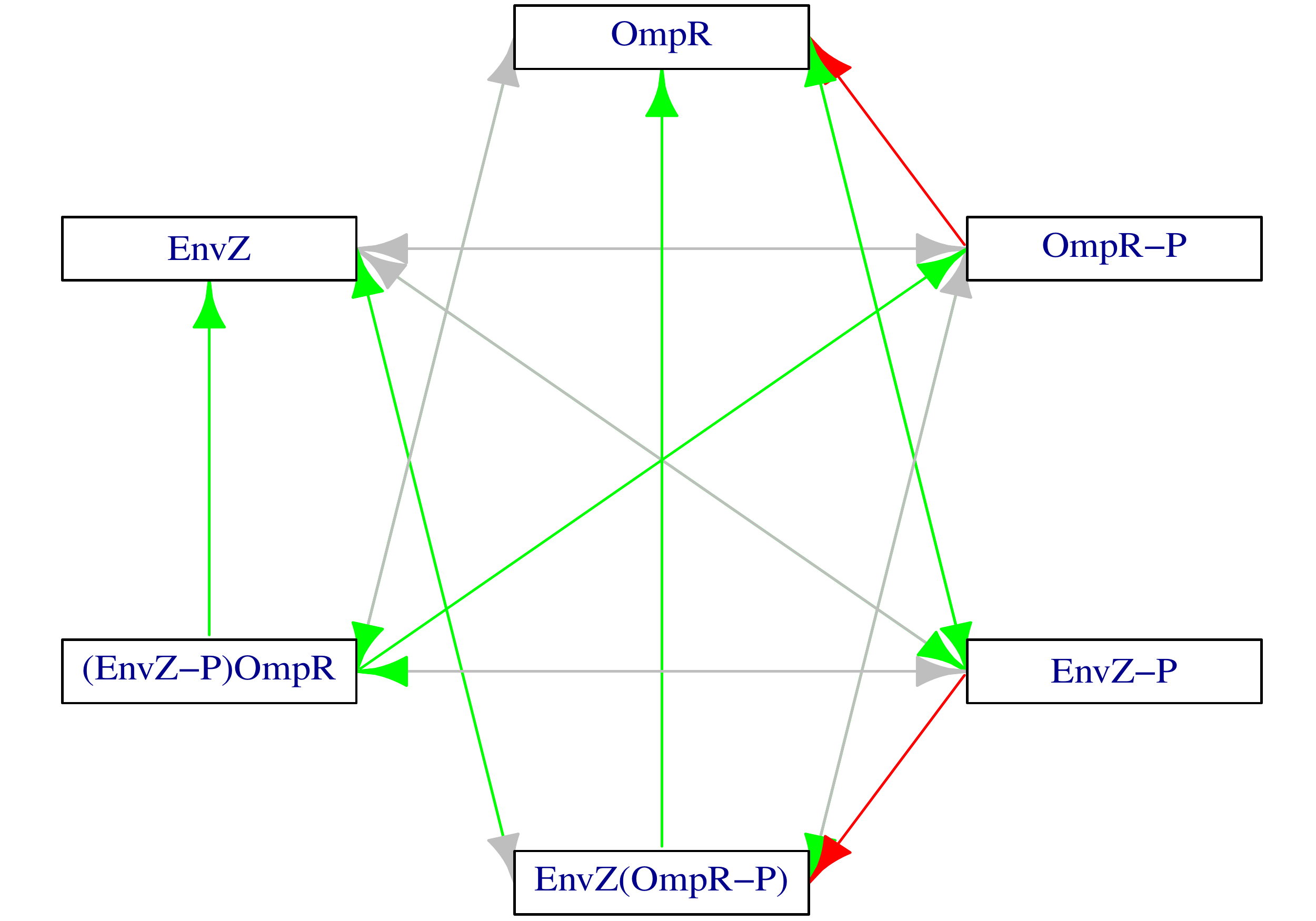}\llap{\parbox[b]{0.3in}{(d)\\\rule{0ex}{1.3in}}} }\\
\end{tabular}
\caption{Simulated time course data (a) from the EnvZ/OmpR system with two perturbations and two interventions.
Networks estimated from the perturbed data (b), the intervened data (c) and all data (d). Edges colouring scheme: true positive (green), false positive (red), false negative (gray).}
\label{fig:EnvZOmpR_graph}
\end{figure}

The EnvZ/OmpR system is driven by the six coupled ordinary differential equations (see e.g., \cite{Batchelor2003}):
\begin{equation}
 \label{eq:Envzequations}
 { \small
 \begin{aligned}
    \frac{d[\textnormal{(EnvZ-P)OmpR}]}{dt}  =& k_1[\textnormal{EnvZ-P}][\textnormal{OmpR}]- (k_{-1}+k_t)[\textnormal{(EnvZ-P)OmpR}] \\
    \frac{d[\textnormal{EnvZ(OmpR-P)}]}{dt}  =& k_2[\textnormal{EnvZ}][\textnormal{OmpR-P}]- (k_{-2}+k_p)[\textnormal{EnvZ(OmpR-P)}] \\
    \frac{d[\textnormal{EnvZ-P}]}{dt}  =& k_{-1}[(\textnormal{EnvZ-P})\textnormal{OmpR}]- k_{-k}[\textnormal{EnvZ-P}] + k_k[\textnormal{EnvZ}] \\
    &-k_1[\textnormal{EnvZ-P}][\textnormal{OmpR}]\\
    \frac{d[\textnormal{EnvZ}]}{dt}  =& k_{-k}[\textnormal{EnvZ-P}] - k_k[\textnormal{EnvZ}]+ (k_p+k_{-2})[\textnormal{EnvZ(OmpR-P)}] \\
    &+k_t[\textnormal{(EnvZ-P)OmpR}]     -    k_2[\textnormal{EnvZ}][\textnormal{OmpR-P}]\\
    \frac{d[\textnormal{OmpR}]}{dt} =&  k_{-1}[\textnormal{(EnvZ-P)OmpR}]-k_{1}[\textnormal{EnvZ-P}][\textnormal{OmpR}] \\
    &+ k_p[\textnormal{(EnvZ)OmpR-P}] \\
    \frac{d[\textnormal{OmpR-P}]}{dt} =&  k_{t}[\textnormal{(EnvZ-P)OmpR}]-k_{2}[\textnormal{EnvZ}][\textnormal{OmpR-P}] \\
    &+ k_{-2}[\textnormal{EnvZ(OmpR-P)}]     
 \end{aligned}
 }
\end{equation}
which is a \textit{mass action kinetics} (MAK) system. See Section \ref{sec:MAK} for details. In these equations, EnvZ-P and OmpR-P denote the phosphorylation of EnvZ and OmpR, respectively. These systems are characterised by a set of reactions, e.g., $$\textnormal{EnvZ(OmpR-P)}\xrightarrow{k_p}\textnormal{OmpR} + \textnormal{EnvZ}.$$ The AIM algorithm works by searching through a large set of candidate reactions.

Figure \ref{fig:EnvZOmpR_graph} shows simulated data at 26 time points and the network recovered from these data via the AIM algorithm (specifically, Algorithm \ref{alg:AIM_wrap} in Section \ref{sec:proposedmeth}). The networks were recovered from a search space consisting of all MAK systems constructed from reactions on the form
\begin{equation}
\begin{aligned}
&X\rightarrow Y, \quad X+Y\rightarrow Z \quad \text{or} \quad Z\rightarrow X+ Y, \\
\text{with } &X,Y,Z \in\left\{ \begin{aligned}
&\textnormal{EnvZ(OmpR-P)}, \textnormal{(EnvZ-P)OmpR}, \\ 
& \textnormal{EnvZ}, \textnormal{EnvZ-P}, \textnormal{OmpR}, \textnormal{OmpR-P} 
\end{aligned}\right\}.  
\end{aligned}
\end{equation}
The true parameter values in \eqref{eq:Envzequations} were drawn at random from a normal distribution with mean 3 and the initial conditions were drawn uniformly at random from the interval $[5,10]$. The AIM algorithm was here tuned to report reaction networks consisting of eight reactions.

This example illustrates that correct recovery of the network of reactions can benefit from combining several types of data sets. It was thus paramount to develop statistical and computational tools for recovering the network from time course data sampled under different perturbations and/or interventions, thus unifying the estimation process and circumventing the need for highly specific and specialised experiments with individual estimation procedures.

\section{Statistical Framework}
\label{sec:stat}
We consider a $d$-dimensional ODE given by: 
\begin{equation}
 \label{eq:ode}
 \frac{dx}{dt} = f(x(t),\theta), \qquad x(0) = x_0
\end{equation}
with initial condition $x_0\in \mathbb{R}^d$ and the smooth field $f:\mathbb{R}^d\times \mathbb{R}^p\rightarrow\mathbb{R}^d$ parameterised by $\theta\in\mathbb{R}^p$. In terms of $f$ we define a corresponding network with nodes $1, \ldots, d$ and an edge from node $l$ to node $i$ if and only if $\partial f_i/\partial x_l \neq 0$. For many parameterised ODE systems a nonzero coordinate in $\theta$ corresponds to the presence of one or a few edges in the network, thus if we enforce sparsity in $\theta$ we also enforce sparsity in the network. This is, for instance, the case for the polynomial and rational fields that are currently implemented in the R package \textit{episode}, see Table \ref{tab:aim}. In the setting of this paper, the focus is therefore on $p$ being large but the true parameter being sparse. In some of the examples we consider, $p$ is of the order $12,000$ with $\theta$ having as little as $0.65\%$ of the parameters being nonzero.

We assume that the process $x$ is observed at discrete time points $(t_i)_{i=1}^n$ with i.i.d. noise $(\varepsilon_i)_i$,
$$y(t_i) = x(t_i) + \varepsilon_i.$$ 
Using a sparsity enforcing penalty function $\mathrm{pen}$, e.g., $\ell^1$, elastic net, SCAD or MCP, we will consider the penalised least squares loss function
\begin{equation}
 \label{eq:loss}
 \ell_y(\theta) \coloneqq \frac{1}{2}\sum_{i=1}^{n}{ \sum_{l=1}^{d} { w_{i,l}(y_l(t_i) - x_l(t_i,\theta))^2 } } +\lambda\sum_{j=1}^{p}{v_j \mathrm{pen}(\theta_j)},
\end{equation}
where $v=(v_j)_j$ are penalty weights and $w=(w_{i,l})_{i,l}$ are observation weights. Strong distributional assumptions on the errors, $\varepsilon_i \in \mathbb{R}^d$, are not necessary, but we note that the least squares loss doesn't account for potential correlation among the $d$ coordinates. However, differences in the error variances among the coordinates are accounted for by the observation weights, which are chosen adaptively by the proposed AIM algorithm.

Finally, we allow for observations of the same system under different interventions. We assume that the interventions are encoded in the ODE system through a Hadamard product of the parameter $\theta$. More precisely, let $\{1,...,E\}$ be a finite set of environments representing the interventions and let the data $y$ consist of $E$ sub-datasets $(y^e)_{e=1}^E$ with $n_e$ time points in environment $e$. Define the environment specific observation weights similarly. The effective parameter of the ODE system in environment $e$ is $\theta \circ c_{e}$, where $\theta\in\mathbb{R}^p$ is the baseline parameter corresponding to the unconstrained/unintervened system and $c_e\in \mathbb{R}^p$ is a vector of coordinatewise scale factors.

Typically, the scale factors $c_e$ are binary. For instance, if in environment $e$ the $l^{\mathrm{th}}$ coordinate of $x$ is inhibited from affecting the $i^{\mathrm{th}}$ coordinate, then coordinate $j$ of $c_e$ is set to zero if and only if $\partial^2f_i/\partial\theta_j\partial x_l \neq 0$. This inhibiting mode-of-action of an intervention is commonly used in gene regulatory networks in which certain proteins can inhibit the translation of some genes (see e.g., \cite{Fire:1999}, \cite{Elbashir:2001}). The loss function taking this type of intervention into account thus reads
\begin{equation}
 \label{eq:loss_intervene}
 \ell_y(\theta) \coloneqq \frac{1}{2}\sum_{e=1}^E\sum_{i=1}^{n_e}{ \sum_{l=1}^{d} { w^e_{i,l}\left(y^e_l(t_i) - x^e_l(t_i,\theta\circ c_e)\right)^2 } } +\lambda\sum_{j=1}^{p}{v_j \mathrm{pen}(\theta_j)}.
\end{equation}
Direct optimisation of \eqref{eq:loss_intervene} is challenging as this is generally a non-convex optimisation problem with many local minima, and most nonlinear ODEs will have to be solved numerically just to evaluate \eqref{eq:loss_intervene}. In the following section we will introduce methods that mitigate some of the difficulties. 

\section{Methods}
\label{sec:meth}
\subsection{The least squares method}
Direct minimisation of \eqref{eq:loss_intervene} above is called the (penalised) \textit{least squares} method. This is sometimes referred to as the \textit{gold standard} approach, see, e.g., \cite{Chen:2016}. As noted above, 
evaluating $x$ in \eqref{eq:loss_intervene} typically requires a numerical ODE solver, which makes the least squares method computationally heavy. We refer to \cite{Sauer:2006} for a comprehensive overview of numerical ODE solvers, and to Appendix \ref{sec:comp_asp} for details on how to optimise \eqref{eq:loss_intervene} while keeping computation time to a minimum.

\subsubsection{Issues}
The penalised least squares method suffers from three main problems: it is computational demanding, it is a non-convex optimisation problem, and the solution depends on the choice of parameter scale (the choice of penalty weights). 

The numerical solution of \eqref{eq:ode} is fundamentally a sequential problem, thus each evaluation of $x$ is computationally heavy with only limited parallelisation options. Moreover, the derivative of $x$ with respect to $\theta$ or $x_0$ solves another ODE, called the \textit{sensitivity equations}, of dimensions $d^2$ and $dp$, respectively (see Appendix \ref{sec:comp_asp} for details).  

The loss function \eqref{eq:loss} is non-convex even in the simplest case of a linear ODE, since linearity of the vector field $f$ does not imply that the solution to the ODE is linear. For nonlinear ODE systems we cannot even expect that \eqref{eq:loss_intervene} has a unique local minimiser for small $\lambda$. 

The dependence on parameter scale is a general problem for penalised nonlinear least squares. The scales on which the parameters are penalised are essential for what parameters the sparsity inducing penalty selects. All other equal, parameters for which $x$ is more sensitive is typically chosen over those for which $x$ is less sensitive. This is a clear issue for correct network inference. In linear regression, it is common to standardise the predictors to bring the parameters on a common scale, but no immediate method exists for standardising the parameters in the nonlinear least squares function \eqref{eq:loss_intervene}. It appears that any such method would depend on the unknown $\theta$. 

The inverse collocation methods introduced below address the three main problems of the least squares method.

\subsection{Inverse collocation methods}
In numerical analysis, \textit{collocation methods} are a class of methods for solving ODE systems numerically. It goes as follows; in the ODE system 
\begin{equation}
 \label{eq:ode_col}
 \frac{dx}{dt} = f(x(t),\theta), \qquad x(0) = x_0
\end{equation}
the parameter vector $\theta$ is assumed known. Moreover, a finite set of collocation points $\mathcal{C} \subseteq \mathbb{R}$ are chosen, as well as a vector space $\mathcal{V}$ of functions. A numerical solution, $\tilde{x} \in \mathcal{V}$, is sought that makes  
$\|\frac{d\tilde{x}}{dt}(t) - f(\tilde{x}(t), \theta)\|$
small in the collocation points for some norm. That is, the numerical solution is found by minimising a distance between 
$\left(\frac{d\tilde{x}}{dt}(t)\right)_{t \in \mathcal{C}}$ and $(f(\tilde{x}(t), \theta))_{t \in \mathcal{C}}$. Typically, $\mathcal{V} = \mathrm{span}(\varphi_j)$ for a choice of finitely many basis functions $\varphi_j:\mathbb{R}\rightarrow \mathbb{R}$ and the norm is the 2-norm on $\mathbb{R}^d$. Collocation methods thus solve the forward problem of computing the solution of \eqref{eq:ode_col} for a known $\theta$. 

By \textit{inverse collocation methods} we refer to a class of estimators of the parameter $\theta$, given the observed trajectory $x$, that solve the inverse problem using the collocation idea. These methods exist in many versions (\cite{Varah:1982}, \cite{Brunel:2008}, \cite{Liang:2008},  \cite{Calderhead:2009}, \cite{Gugushvili:2012}, \cite{Dondelinger:2013}) and are known under many other names, e.g., \textit{gradient matching}, \textit{trajectory matching} or \textit{smooth-and-match estimators}. However, they all rely on the same two-step procedure: 1) approximate the data, $y$, by an element in $\mathcal{V}$ to get an estimate of the full trajectory $\hat{x}$; 2) base the estimation of $\theta$ on the trajectory $\hat{x}$ as if it were the true trajectory, by minimising the distance between the position, gradient or integral at a given set of collocation points. Typically, $\hat{x}$ is obtained as a smoother or an approximation of $y$ via a basis expansion.

One example of an inverse collocation method is the \textit{gradient matching method} (\cite{Varah:1982}, \cite{Brunel:2008}), which minimises the approximate loss function:
\begin{equation}
 \label{eq:approx_loss_gradient}
 \frac{1}{2}\sum_{t \in \mathcal{C}}{\left\| \frac{d\hat{x}}{dt}(t) - f(\hat{x}(t), \theta) \right\|_2^2}.
\end{equation}

This method considerably reduces the computational cost compared to the least squares method, since it does not require solving the ODE system. Moreover, if $f$ is linear in $\theta$ the optimisation problem becomes a linear least squares problem, which thus avoids all the three problems with the least squares method.

Later \cite{Dattner:2015} proposed minimising  
\begin{equation}
 \label{eq:approx_loss_integral}
 \frac{1}{2}\sum_{t \in \mathcal{C}}{\left\| \hat{x}(t) - x_0 - \int_{0}^t{f(\hat{x}(s), \theta) \ ds} \right\|_2^2},
\end{equation}
since the ODE system can be characterised as solving 
\begin{equation}
 \label{eq:inte}
 x(t_2) - x(t_1) = \int_{t_1}^{t_2}{f(x(s),\theta) \ ds}, \qquad \text{for all} \ t_1,t_2\in\mathbb{R}.
\end{equation}
instead. This requires numerical integration, which is often more stable than numerical differentiation, and under certain assumptions $\sqrt{n}$-consistency is guaranteed, as by \cite{Gugushvili:2012}. Also, in this method the smoothed trajectory $\hat{x}$ does not have to be differentiable.

Note that in all of the above methods the collocation time points in $\mathcal{C}$ do not have to coincide with the observation time points of $y$. However, adding more time points in \eqref{eq:approx_loss_gradient} and \eqref{eq:approx_loss_integral} will not necessarily decrease the variance of the estimator, as that mostly comes down to the $y$-$\hat{x}$ relation, i.e., the smoothing operation.

Nonparametric inverse collocation methods also exist, most notably are those by \cite{Wu:2014} and \cite{Chen:2016}. Here the authors do not assume a parametric form of the field $f$, but approximate it by a nonparametric basis. In the former the authors consider the loss function 
\begin{equation}
 \label{eq:approx_loss_gradient_np}
 \frac{1}{2}\sum_{t \in \mathcal{C}}\sum_{l=1}^d {\left(\frac{d\hat{x}_l}{dt}(t) - \sum_{j,k}\psi_k(\hat{x}_j(t))\theta_{ljk}\right)^2},
\end{equation}
with $(\psi_k)_{k=1}^K$ a finite set of basis functions and $(\theta_{ljk})_{ljk}$ estimable parameters. In \cite{Chen:2016} the integrals are considered instead:
\begin{equation}
 \label{eq:approx_loss_integral_np}
 \frac{1}{2}\sum_{t \in \mathcal{C}}\sum_{l=1}^d {\left(\hat{x}_l(t) - x_l(0)- \sum_{j,k}\Psi_k(\hat{x}_j)(t)\theta_{ljk}\right)^2},
\end{equation}
with $\Psi_k(x)(t)\coloneqq \int_{0}^{t}{\psi_k(x(s)) \ ds}$. Note that both nonparametric methods assume $f$ to be additive in the coordinates of $x$.

Finally, we note that the generalised profiling method described by \cite{Ramsay:2007} is another variation on the inverse collocation method. It is inspired by functional data analysis and the main difference lies in that the smoothing step is $\theta$-dependent and thus becomes part of the optimisation step.

Penalised versions of the inverse collocation methods -- as alternatives to minimising \eqref{eq:loss} -- have also been proposed by e.g., \cite{Lu:2011} and \cite{Wu:2014} to promote sparse solutions.

\subsubsection{Issues}
Though the inverse collocation methods remedy most issues of the least squares approach (in fact all of those discussed above, if the ODE is $\theta$-linear), the inverse collocation methods also have their share of issues. Most notably, the results become dependent on the initial approximation (the smoother), which will introduce a bias without a clear trade-off in terms of a reduced variance. To illustrate this we present a small simulation study. Consider the classic Michaelis-Menten kinetics modelling the chemical reaction system (see \cite{MM:1913})
\begin{equation}
 \label{eq:MM_reaction}
 \mathrm{E} + \mathrm{S} \xrightleftharpoons[k_r]{k_f} \mathrm{ES} \xrightarrow{k_{cat}} \mathrm{E} + \mathrm{P}
\end{equation}
in which the enzyme (E) forms a complex (ES) through a binding interaction with the substrate (S), which further releases the product (P) along with the freed enzyme. The abundances of the four compounds $x=(x_{\mathrm{E}}, x_{\mathrm{ES}}, x_{\mathrm{P}}, x_{\mathrm{S}})$ satisfy an ODE with $p=3$ positive parameters $(k_f, k_r, k_{cat})$:
\begin{equation}
\label{eq:MM_ODE}
  \begin{aligned}
    \frac{dx_{\mathrm{E}}}{dt} = -k_fx_{\mathrm{E}}x_{\mathrm{S}} + k_rx_{\mathrm{ES}} + k_{cat}x_{\mathrm{ES}} & \qquad \frac{dx_{\mathrm{P}}}{dt} = k_{cat}x_{\mathrm{ES}} \\
    \frac{dx_{\mathrm{ES}}}{dt} = k_fx_{\mathrm{E}}x_{\mathrm{S}} - k_rx_{\mathrm{ES}} - k_{cat}x_{\mathrm{ES}} & \qquad    \frac{dx_{\mathrm{S}}}{dt} = -k_{f}x_{\mathrm{E}}x_{\mathrm{S}} + k_r x_{\mathrm{ES}}. 
  \end{aligned}
\end{equation}

\begin{figure}
  \centering
    \includegraphics[width=0.95\textwidth]{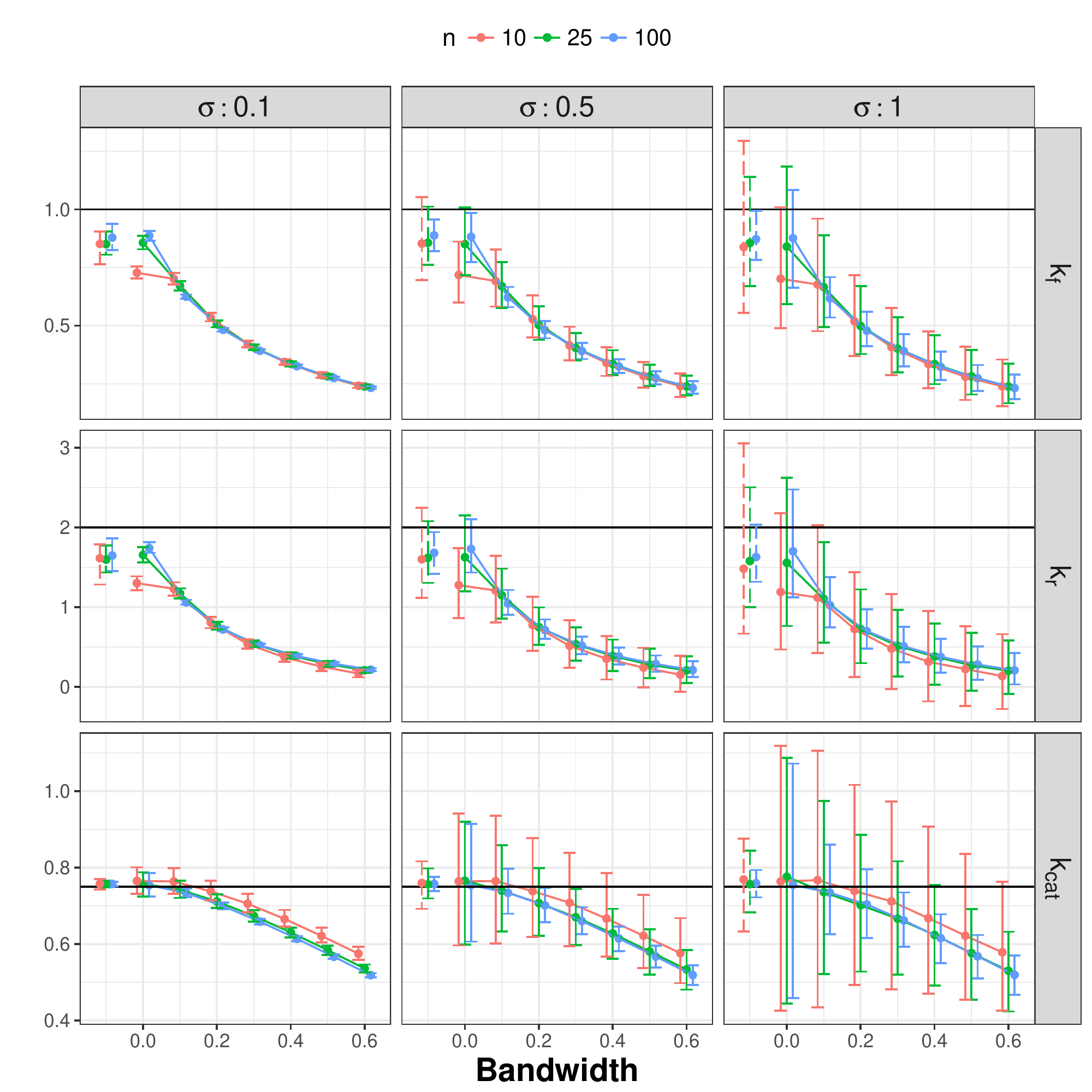}
  \caption{Medians and $5\%$ and $95\%$ percentiles of inverse collocation estimator considered by \cite{Dattner:2015}. The kernels were scaled such that the quartiles are at $\pm 0.25 \times \mathrm{bandwidth}$. Data is simulated from \eqref{eq:MM_ODE} with $x_0 = (10, 2, 2, 10)$ and a time range of 1. The true parameters are marked with horizontal lines. The dashed lines on the left are the corresponding medians and percentiles of the least squares method.}  
  \label{fig:MM}
\end{figure}

This classical ODE model is linear in the parameters and thus well suited for the inverse collocation methods. We generated data at $n=10, 25, 100$ equidistant time points from the true trajectory with i.i.d. additive Gaussian noise. The data set was replicated $250$ times and for each of them we applied a Gaussian kernel smoother with a range of bandwidths followed by the method proposed by \cite{Dattner:2015} to obtain parameter estimates. A summary of the resulting estimators is presented in Figure \ref{fig:MM}. 

From Figure \ref{fig:MM} we notice a bias which severely increases with the bandwidth, while the variance is only moderately reduced. Moreover, the bias hardly seems to change with the number of observations, unless the bandwidth is zero (equivalent to a linear interpolation smoother). Intuitively, this is no surprise: the purpose of smoothers, as indicated by their name, is to smooth the data. This is often manifested in a reduced pointwise variance, $\mathrm{V}(\hat{x}_{y}(t))\leq \mathrm{V}(y(t))$ for $t\in\mathbb{R}$, and an increased autocovariance, $\mathrm{Cov}(\hat{x}_{y}(t), \hat{x}_{y}(s))\geq \mathrm{Cov}(y(t), y(s))$, for $t,s\in\mathbb{R}$ close. Together this results in underestimated slopes. Since the slopes are essentially what is being modelled in ODE systems we would expect the resulting parameter estimates to have a large bias. 

The least squares method and inverse collocation with zero bandwidth have the smallest biases. However for moderate and large noise levels the variance of the least squares method decreases faster with the number of observations. Though the inverse collocation methods with large bandwidths have slightly smaller variance, the least squares method still outperforms them, except for some settings with $n=10$ and $\sigma=0.1$. 

Finally, inverse collocation methods suffer from one additional issue; they require fully observed processes to work. There is no obvious way of producing smoothed curves for latent coordinates and all coordinates are required in \eqref{eq:approx_loss_gradient} and \eqref{eq:approx_loss_integral}. This problem is revisited in Section \ref{sec:glyco}.


\subsection{Adaptive Integral Matching}
\label{sec:proposedmeth}
We propose combining an inverse collocation method with the least squares method in such a way that we benefit from both methods. Inverse collocation methods are computationally lighter and produce good approximate parameter estimates, while not fully enjoying the statistical qualities of the least squares estimator. The least squares method is computationally expensive and suffers heavily from multiple local minimas, while generally performing better if the latter problems are alleviated. 

Before presenting our suggestion of a combined estimator, we introduce a modification of the inverse collocation method by \cite{Dattner:2015}. We propose the collocation method that consists of minimising the following approximate loss function
\begin{equation}
 \label{eq:approx_loss_aim}
 \begin{aligned}
 \tilde{\ell}^{\hat{x}}(\theta) \coloneqq & \frac{1}{2}\sum_{e=1}^E\sum_{i=1}^{n_e-1}{ \sum_{l=1}^{d}{ w^e_{i,l}\left( \hat{x}^e_l(t_{i+1}) - \hat{x}^e_l(t_{i}) - \int_{t_i}^{t_{i+1}}{f_l(\hat{x}^e(s), \theta \circ c_e) \ ds}\right)^2}} \\
 & +\lambda\sum_{j=1}^{p}{v_j \mathrm{pen}(\theta_j)}, 
 \end{aligned}
\end{equation}
where $\hat{x}^e$ is the smoothed curve based on the data from environment $e$, and $\hat{x}=(\hat{x}^e)_{e=1}^E$ denotes the collection of smoothed curves for each environment. The above differs from \eqref{eq:approx_loss_integral} by integrating between consecutive time points instead of integrating from $0$ to $t$. This has two positive side effects: 1) it prevents errors between the true trajectory and its estimate $\hat{x}$ from accumulating; 2) the initial condition $x_0$ is no longer estimated. This is highly preferable as the initial condition is often a nuisance parameter and in a penalised setup the optimisation procedure often sets $x_0$ to compensate for the restricted freedom in $\theta$. We refer to the estimator
\begin{equation}
 \label{eq:int_match}
  \hat{\theta}^{\hat{x}}_{\lambda} \coloneqq \argmin_{\theta} \tilde{\ell}^{\hat{x}}(\theta)
\end{equation}
as the \textit{integral matching estimator} and stress that it depends on the smoother, $\hat{x}$.

%

From an integral matching estimate, $\hat{\theta}^{\hat{x}}_{\lambda}$, we adapt the scales $(c_e)_e$ and, optionally, the weights $(w^e)_e$. The new adapted scales are proportional to 
\begin{equation}
\label{eq:adapt_sc} 
  c_e \circ \left(\left\| \left( \int_{t_i}^{t_{i+1}}{ \partial_{\theta_j}f(\hat{x}^e(s), \hat{\theta} \circ c_e) \ ds} \right)_{i,e} \right\|^{-1}_2\right)_j, \quad \text{for } e=1,...,E.   
\end{equation}
If the field is linear in $\theta$, then the updated scales only depend on the smoother. If the field is not $\theta$-linear one uses $\hat{\theta}=\hat{\theta}^{\hat{x}}_{\lambda}$ for a small $\lambda$. The scales are thus standardised by the column norms of the first order Taylor approximation of the integrals in \eqref{eq:approx_loss_aim}. If $f$ is linear in $\theta$, this coincides with standardising the columns in a penalised linear least squares problem. This adaptation of the scales ensures that parameters are locally on the same scale and thus penalised in a fair manner in the subsequent least squares estimation. Similarly, the new adapted weights are proportional to
\begin{equation}
\label{eq:adapt_weight}
  \frac{(w_{i,l}^e)_{i,e}}{\sum_{e=1}^E\sum_{i=1}^{n_e-1}{ { w^e_{i,l}\left( \hat{x}^e_l(t_{i+1}) - \hat{x}^e_l(t_{i}) - \int_{t_i}^{t_{i+1}}{f_l(\hat{x}^e(s), \theta \circ c_e) \ ds}\right)^2}}}
\end{equation}
for $l=1,...,d$, i.e., inversely proportional to the empirical variances for each species. This adapts the variance structure across species for the subsequent estimation. This leads to the \textit{adaptive integral matching} (AIM) algorithm:
\begin{alg}[AIM]
\label{alg:AIM} 
Input: Time course data from $E$ environments, $y=(y^e)_{e=1}^E$, each sampled at $(t_i)_{i=1}^{n_e}$ timepoints. Similarly structured observation weights $w=(w^e)_{e=1}^E$, along with penalty weights, $v\in \mathbb{R}^p_+$, and environment-specific scales $(c_e)_{e=1}^E$. Smoothed trajectories $(\hat{x}^e)_e$ evaluated on a fine grid of time points.
\begin{enumerate}
 \item \label{itm:im} Apply the integral matching estimator, \eqref{eq:int_match}, to obtain initial estimates $\hat{\theta}^{\hat{x}}_{\lambda}$ for a sequence of $\lambda$ values.
 \item Adapt the scales and weights according to \eqref{eq:adapt_sc} and \eqref{eq:adapt_weight}.
 \item \label{itm:refit} Refit by minimising \eqref{eq:loss_intervene} using the adapted weights and $\hat{\theta}^{\hat{x}}_{\lambda}$ as initial value.
\end{enumerate} 
\end{alg}

In step \eqref{itm:refit} the penalty term may be scaled down or removed entirely to reduce the bias induced by the penalty, and the parameter space may be restricted to reduce the computational costs. Algorithm \ref{alg:AIM_wrap} below presents a particular incarnation of the refitting step. In Appendix \ref{sec:comp_asp} additional techniques to reduce the computation time are presented.

%
%

\subsubsection{Implementation}
\label{sec:implementation}
As part of this paper, software for optimising \eqref{eq:loss_intervene} and \eqref{eq:approx_loss_aim} (used in Algorithm \ref{alg:AIM}) is available in the R package \textit{episode}. In the latter optimisation problem the user supplies the smoothed trajectories $(\hat{x}^e)_e$ evaluated on a fine grid of time points and the software then optimises \eqref{eq:approx_loss_aim} using numerical integration over the supplied grid. By keeping this modular form, the user has complete freedom in choosing a suitable smoother. It is possible not to smooth the data at all, which corresponds to $\hat{x}^e$ linearly interpolating the observations $y^e$. 

\begin{figure}
  \centering
    \includegraphics[width=0.9\textwidth]{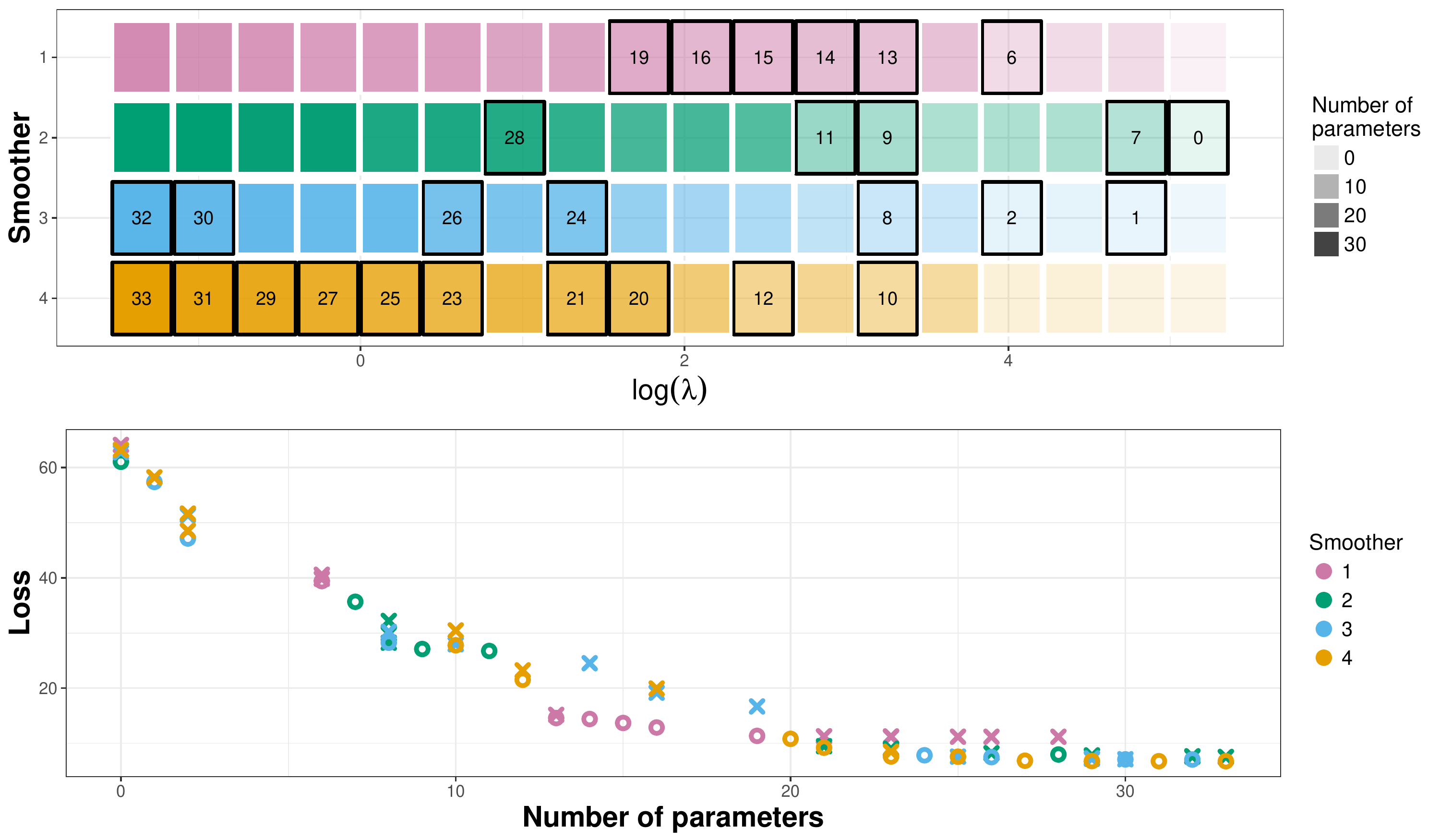}
  \caption{Visualisation of the stratified ranking in Algorithm \ref{alg:AIM_wrap} applied to the EnvZ/OmpR data from Section \ref{sec:netw}. Four smoothers were employed and each proposed a sequence of candidate models for varying tuning parameter (top). The loss values of the candidate models are stratified according to model size (bottom). For each model size the candidate with minimal loss is found and marked with a black border (top).}  
  \label{fig:Alg_visual}
\end{figure}

We recommend subjecting a whole family of smoothers to Algorithm \ref{alg:AIM} in order to alleviate potentially high variance and multiple local minima issues. The resulting version of the AIM algorithm that we suggest and have tested extensively consists of the following steps:

\begin{alg}
\label{alg:AIM_wrap} 
Input: Time course data from $E$ environments, $y=(y^e)_{e=1}^E$, each sampled at $(t_i)_{i=1}^{n_e}$ timepoints. Similarly structured observation weights $w=(w^e)_{e=1}^E$, along with penalty weights, $v\in \mathbb{R}^p_+$, and environment-specific scales $(c_e)_{e=1}^E$. 
\begin{enumerate}
 \item Produce a family of smoothed curves $\{\hat{x}\}$, from data $y$, where the smoother is applied to each environment separately: $\hat{x}=(\hat{x}_e)_{e=1}^E$.
 \item For each $\hat{x}$ apply Algorithm \ref{alg:AIM} with the refitting step implemented as follows: define the support estimator $\hat{S}^{\hat{x}}_{\lambda}=\mathrm{supp}(\hat{\theta}^{\hat{x}}_{\lambda})$ and compute the unpenalised least squares estimate
 \begin{equation}
  \label{eq:loss_no_pen}
  \tilde{\theta}^{\hat{x}}_{\lambda} \coloneqq \underset{\theta:\mathrm{supp}(\theta)=\hat{S}^{\hat{x}}_{\lambda}}{\arg\min}\frac{1}{2}\sum_{e=1}^E\sum_{i=1}^{n_e}{ \sum_{l=1}^{d} { w^e_{i,l}\left(y^e_l(t_i) - x^e_l(t_i,\theta\circ c_{e})\right)^2 } }.
 \end{equation}
over the restricted parameter space determined by $\lambda$ and $\hat{x}$. 
 \item Stratify the refitted estimates $(\tilde{\theta}^{\hat{x}}_{\lambda})_{\lambda, \hat{x}}$ by the number of non-zero parameters. For each strata rank the resulting estimates by their loss value at optimum. See Figure \ref{fig:Alg_visual} for an illustration of this step. 
\end{enumerate} 
\end{alg}

The purpose of the stratified ranking is to produce a sequence of models indexed by the number of nonzero parameters. This is primarily important for comparison purposes in the subsequent sections. 

Currently, the R package \textit{episode} implements AIM and other learning algorithms for mass action kinetics (described below), which encode all polynomial fields, power law kinetics, which encode all polynomial fields in a different way and two larger classes of ODE systems assuming a rational form of the field. As for penalties, $\ell^1$, $\ell^2$, elastic net, SCAD and MCP are implemented. Moreover, the package handles missing values and allows for box-constrained optimisation as well. Table \ref{tab:aim} provides a schematic overview of the features in \textit{episode}. The tools in \textit{episode} are flexible and modular and the Algorithms \ref{alg:AIM} and \ref{alg:AIM_wrap} are primarily recommendations on how to combine them. When using the \textit{episode} package for the least squares method, i.e., optimising \eqref{eq:loss_intervene}, suitable initialisations are required and the resulting estimates may depend on these. The tools are thus designed to easily pass the  integral matching estimates as initialisations for the least squares method.

{
\scriptsize
\begin{table}
  \taburowcolors[2] 2{tableLineOne .. tableLineTwo}
  \tabulinesep = ^2mm_1mm
  \everyrow{\tabucline[.2mm  white]{}}  
  \begin{tabu} to \textwidth {l >{\bfseries}X[r, 1] X[3] l}
    \tableHeaderStyle
    &  & ODE Models & \\
    & MAK & Mass Action Kinetics 
    $$
    \frac{dx}{dt} = (B-A)^T\mathrm{diag}(x^A)k,
    $$
    with $A, B\in\mathbb{N}_0^{r\times d}$ fixed and $k\in\mathbb{R}^r_+$ estimable.
    & \\
    & PLK & Power Law Kinetics 
    $$
    \frac{dx}{dt} = \theta x^A,
    $$
    with $A\in\mathbb{N}_0^{r\times d}$ fixed and $\theta\in\mathbb{R}^{d\times r}$ estimable.
    & \\
    & RLK & Rational Law Kinetics 
    $$
    \frac{dx}{dt} = \theta \frac{x^A}{1+x^B},
    $$
    with $A, B\in\mathbb{N}_0^{r\times d}$ fixed, the fraction  evaluated elementwise and $\theta\in\mathbb{R}^{d\times r}$ estimable.
    & \\
    & RMAK & Rational Mass Action Kinetics 
    $$
    \frac{dx}{dt} = C^T \frac{\theta_1x^A}{1+\theta^2x^A},
    $$    
    with $A\in\mathbb{N}_0^{b\times d}$ and $C\in\mathbb{N}_0^{r\times d}$ fixed, the fraction evaluated elementwise and $\theta_1,\theta_2\in\mathbb{R}^{r\times b}$ estimable.
    & 
  \end{tabu} \vspace{1mm}
  
  \taburowcolors[2] 2{tableLineOne .. tableLineTwo}
  \tabulinesep = ^2mm_1mm
  \everyrow{\tabucline[.2mm  white]{}}
  \begin{tabu} to \textwidth {l >{\bfseries}X[r, 1] X[3] l}
    \tableHeaderStyle
    &  & Data Structures  & \\
    & Inhibition & Species $i$ is inhibited from reacting with species $j$ in environment $e$: Set $l^{\mathrm{th}}$ coordinate of $c_e\in\mathbb{R}^p$ to $0$ if $\partial_{\theta_l}f_{ij} \neq 0$, and 1 otherwise.
    & \\
    & Activation & Species $i$ only reacts with species $j$ in environment $e$: If $\partial_{\theta_l}f_{ij}\neq 0$ set $c_e(l)=1$ and $c_{e'}(l)=0$ for all $e\neq e'$.
    & \\
    & Stimulation & Reaction rate of reaction $l$ is increased by factor $k$ in environment $e$: Set $c_e(l)=k$. 
    & \\
    & Misc & Missing data. Partially observed processes only supported by exact estimation.
    & 
  \end{tabu} \vspace{1mm}
  
  \taburowcolors[2] 2{tableLineOne .. tableLineTwo}
  \tabulinesep = ^2mm_1mm
  \everyrow{\tabucline[.2mm  white]{}}
  \begin{tabu} to \textwidth {l >{\bfseries}X[r, 1] X[3] l}
    \tableHeaderStyle
    &  & Estimation Components & \\
    & Penalties & $\ell^1$, $\ell^2$, elastic net, SCAD, MCP and no penalty.
    & \\
    & Weights & Both observation and penalty weights.
    & \\
    & Loss & Can minimise both least squares loss \eqref{eq:loss_intervene} and integral matching loss \eqref{eq:approx_loss_aim}. The minimiser of the latter loss function can easily be passed as initialisation for minimising the former.
    & \\
    & Parameter Constraints & Box constraints for all estimable parameters are available.
    & \\
    & Misc. & Automatic adaptation of parameter scales and observation weights. 
    & 
  \end{tabu} 
\caption{Overview of features in the R package \textit{episode}.}
\label{tab:aim}
\end{table}
}
 
\clearpage

\section{Applications}
\label{sec:app}
In this section we study two concrete large scale dynamical systems. One is the \textit{in silico} protein phosphorylation network used in the eighth DREAM challenge, and the other is glycolysis in \textit{Saccharomyces cerevisiae}. Both of these systems are like the EnvZ/OmpR system based on \textit{mass action kinetics}, which is first reviewed briefly. However, in these applications it is not all components of the mass action system that is observed, and rational fields are used to model the dynamics of the observed species. 

\subsection{Mass Action Kinetics}
\label{sec:MAK}
We consider a chemical kinetics framework of ODE systems. Assuming that we have $d$ chemical species, e.g., $\mathrm{NaCl}$, $\mathrm{H}_2\mathrm{O}$ or proteins, labelled $X=(X_1,...,X_d)$. A set of $r=1,...,R$ reactions on the form:
\begin{equation}
 \label{eq:react}
 a_1X_1 + ... + a_dX_d \rightarrow  b_1X_1 + ... + b_dX_d,
\end{equation}
govern the dynamics of the species. Here $(a_i)_{i=1}^d$ and $(b_i)_{i=1}^d$ are non-negative integers, called the stoichiometric coefficients. For reaction $r=1,...,R$ let $A_r, B_r \in \mathbb{N}_0^d$ denote the vector of left hand and right hand side stoichiometric coefficients, respectively. The net change of molecules due to reaction $r$ is $v_r = B_r - A_r$.

Let $x=(x_1,...,x_d)\in\mathbb{R}^d_+$ denote the vector of abundances of each chemical species. If the total number of molecules is sufficiently large, we can model the dynamics of $x$ as
\begin{equation}
 \label{eq:chem_kinetics}
 \frac{dx}{dt} = \sum_{r=1}^R{v_r\gamma_r(x(t))}, \qquad x(0) = x_0.
\end{equation}
See \cite{Wallace:2011} for details on its derivation. The laws of mass action kinetics (see, e.g., \cite{Horn:1972}) states that 
\begin{equation}
 \label{eq:mak_rate}
 \gamma_r(x) = k_r x^{A_r},
\end{equation}
where $k_r\geq 0$ is a rate constant and $x^a$ is shorthand for $\Pi_{i=1}^dx_i^{a_i}$ for any two non-negative vectors in $\mathbb{R}^d$. The \textit{stoichiometric matrices} $A$ and $B$ are the $R\times d$-dimensional matrices with the $r^{\mathrm{th}}$ row being $A_r$ and $B_r$ respectively. The matrix notation of \eqref{eq:chem_kinetics} is
\begin{equation}
 \label{eq:mak_matrix}
 \frac{dx}{dt} = (B-A)^T\mathrm{diag}(x^A)k, \qquad x(0) = x_0,
\end{equation}
where $k=(k_r)_{r=1}^R$ and $x^A=(x^{A_r})_{r=1}^R$. 

Ideally, all chemical reaction systems should approximately be a mass action kinetics system. However, in complex reaction networks this may not be the case for the observable species. For gene regulatory networks, say, some proteins may exist in different forms depending on whether an inhibitor or activator is bound to its associated sites, which is not directly observable. In such cases a \textit{quasi-stationary approximation} is often employed to reduce a full mass action system to a system for the observable variables only. The quasi-stationary approximation assumes that the chemical species, $X$, can be divided into two subsets, $X_L$ and $X_O$, the \textit{latent} and \textit{observed} species:
\begin{equation}
 \label{eq:qstat}
 \begin{aligned}
 \frac{dx_L}{dt} &= f_L(x_L, x_O) \\ 
 \frac{dx_O}{dt} &= f_O(x_L, x_O).
 \end{aligned}
\end{equation}
Under the quasi-stationary assumption, i.e., the dynamics of $x_L$ is faster than $x_O$, the dynamics of $x$ can be approximated by the ODE system
\begin{equation}
 \label{eq:qstat_approx}
 \begin{aligned}
 \frac{dx_O}{dt} &= f_O(\tilde{x}_L(x_O), x_O),
 \end{aligned}
\end{equation}
where $\tilde{x}_L(x_O)$ is the restriction of $x_L$ to the manifold $\mathcal{M}_{x_O}\coloneqq\{x_L \mid f_L(x_L, x_O) = 0 \}$ for all values of $x_O$. In certain settings, including fast binding on gene-sites in gene regulatory networks, this approximation is reasonable and the right hand side of \eqref{eq:qstat_approx} is rational. See \cite{Santillan:2008} for a detailed treatment. This is the main motivation for including rational systems in our framework and in the R package \emph{episode}, and its usage will be illustrated by the two applications below.

\subsection{\textit{in silico} phosphoprotein abundance data}
In this section we compare AIM to state-of-the-art network inference methods in systems biology. The eighth DREAM challenge (\cite{DREAM8}) aimed at advancing causal inference of signalling networks in protein phosphorylation. One of the challenges presented the participants with time course data from a complex \textit{in silico} dynamical model of a protein signalling network. The species were given anonymous labels and thus no prior knowledge of the network was given. 

The data consisted of 20 environments produced using combinations of three inhibitors (or no inhibitor) and two types of stimuli each with two strengths. The targets of the inhibitors were provided and encoded in AIM through the scales $(c_{e})_e$ in Algorithm \ref{alg:AIM_wrap}. In light of the rational ODE systems discussed in Section \ref{sec:MAK}, AIM fitted the ODE system given by the field
\begin{equation}
 \label{eq:rational_mak_restrict1}
 \frac{dx}{dt} = \theta \mathrm{diag}(x^A)\mathrm{diag}(1 + x^B)^{-1},
\end{equation}
where $A$ and $B$ are $R\times d$-dimensional matrices and $\theta\in \mathbb{R}^{d\times R}$ estimable coefficients. The rows of $A$ and $B$ ($(a_r)_r$, $(b_r)_r$), ran over all non-negative integer $d$-tuples summing to at most one. Thus the search space consisted of first order rational functions. 

Besides the final DREAM challenge submissions, AIM was compared to two additional methods. The first was the integral matching (IM) estimator, given in \eqref{eq:int_match}. This method represents the use of a penalised inverse collocation method to select the network. The second method was the least squares estimator using a SCAD penalty (SCAD), which was obtained by optimising \eqref{eq:loss_intervene} for a decreasing sequence of $\lambda$, initialised in $\theta=0$. The continuation principle was used, i.e., the optimum found at the previous value of $\lambda$ was re-used as initialisation for next value of $\lambda$.

\begin{figure}[h]
  \centering
    \includegraphics[width=0.9\textwidth]{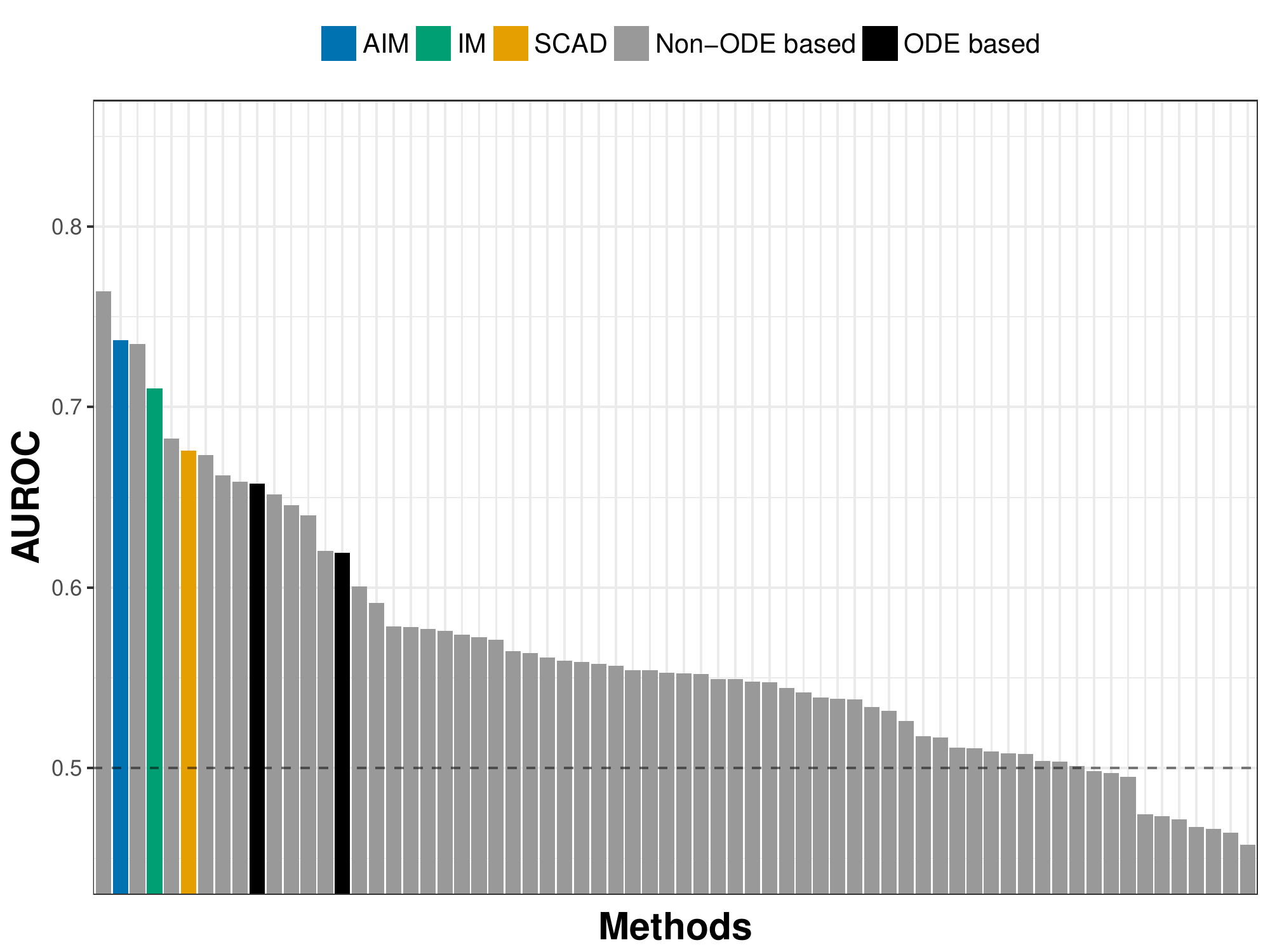}
  \caption{The AUROC scores of all final submissions in the \textit{in silico} network recognition challenge in the eighth DREAM challenge (\cite{DREAM8}) (gray and black bars), along with three of the methods considered in this paper.}  
  \label{fig:dream_aurocs}
\end{figure}

The performance of AIM was assessed using the \textit{DREAMTools} Python package provided by \cite{DREAMtools} and containing the tools used to assess the original challenge submissions. AIM got a AUROC score of $0.737$, which makes AIM the second best solution overall among the 65 submissions and notably better than the two ODE-based submissions. An overview of the performances of AIM, IM and SCAD, along with the final submissions for the eighth DREAM challenge is presented in Figure \ref{fig:dream_aurocs}.

\subsection{Glycolysis in \textit{Saccharomyces cerevisiae}}
\label{sec:glyco}
\cite{Hynne:2001} presented a full scale chemical kinetics model for glycolysis in \textit{Saccharomyces cerevisiae}, constructed from experimental substrate measurements. While \cite{Hynne:2001} knew the metabolic pathway a priori and focused on estimating unknown rate parameters, we will apply AIM to identify the network from simulated data. In total, $d=22$ chemical species enter the glycolysis cycle in an elaborate metabolic network, see Figure \ref{fig:netw_summary}. 

The dynamical model considered by \cite{Hynne:2001} does not fall into the class of mass action kinetics models. All mass action kinetics models have polynomial fields, but the ODE field consider by \cite{Hynne:2001} is rational. More precisely, the field is \eqref{eq:chem_kinetics} with rate functions on the form
\begin{equation}
 \label{eq:rational_mak}
 \gamma_r(x) = \frac{\langle a_r; x^{A_r}\rangle}{\langle b_r; x^{B_r}\rangle}, \qquad A_r\in\mathbb{N}_0^{\alpha_r\times d}, B_r\in\mathbb{N}_0^{\beta_r\times d}, \alpha_r \text{ and } \beta_r \in \mathbb{N},
\end{equation}
with $\langle\cdot;\cdot\rangle$ denoting the standard inner product and $a_r\in \mathbb{R}^{\alpha_r}, b_r\in \mathbb{R}^{\beta_r}$ estimable coefficients. 

For a parametric model on the form \eqref{eq:rational_mak} to be generic enough to include the model considered by \cite{Hynne:2001}, the polynomials $\langle a_r; x^{A_r}\rangle,\langle b_r; x^{B_r}\rangle$ need to have an order of at least 3. Hence, if no prior knowledge on the glycolysis is given, at least $p = 2d(1+d+d^2+d^3) = 490,820$ parameters are needed. It is possible to use AIM with half a million parameters for polynomial systems as given by \eqref{eq:mak_rate}. However, the rational ODE systems are far more sensitive than the polynomial, which in practice results in far longer computations for the numerical solvers and a higher variance of the resulting estimator. Thus a model search space of dimension 490,820 is currently not feasible for rational systems, and we will therefore consider three scenarios for fitting this system using either prior knowledge or an approximate and smaller model search space.

We consider two different prior knowledge scenarios: 1) knowing what \textit{complexes} can be formed in the system, i.e., what terms $A_r\in\mathbb{N}_0^{\alpha_r\times d}, B_r\in\mathbb{N}_0^{\beta_r\times d}$ may appear in the rational field. Even with this prior knowledge, we know very little about the network, since we do not know what complexes drive what reactions. In the system of \cite{Hynne:2001} there are in total 46 complexes. 2) we know a superset of the complexes. In this setting we include an additional 46 false complexes drawn at random. 

In the third scenario we restrict AIM to a smaller parametric model, which will not include the true model. Hence the purpose of this scenario is partly to study the performance of AIM on large and realistic ODE systems and partly to study the robustness to model misspecification. The restricted model space assumes rate functions on the form
\begin{equation}
 \label{eq:rational_mak_restrict}
 \gamma_r(x) = \frac{k_{r}x^{a_r}}{1 + x^{b_r}}, 
\end{equation}
with $a_r\in \mathbb{N}_0^d$ and $b_r\in \mathbb{N}_0^d$ covering all first order terms (i.e., all combinations of non-negative integers summing to at most $1$) and $k_r$ estimable coefficient. This produces a total of $p=d(d+1)^2=11,638$ parameters. By assuming fixed coefficients in the denominator of the rate functions, we obtain an ODE field which is linear in the parameters.

\subsubsection{Simulation study design}
Using the reactions and rate functions listed in Table 1 and 2 in \cite{Hynne:2001}, we numerically solved the ODE system with parameters in Tables 4-7 in \cite{Hynne:2001}. 

We considered $E = 5, 10, 15, 20$ environments each given its own inhibition. These were produced as follows: $20$ distinct chemical species were selected at random, one for each of the maximal number of environments. In each environment the selected species were inhibited, i.e., the species did not form any complexes with the other species and were thus prevented from reacting with the other species.

The trajectories ran for $5$ minutes, at which the system had settled at a stationary point. The trajectories were observed at $30$ log-equidistant time points with additive Gaussian noise, with standard deviations $\sigma = 0.1, 0.25, 0.5$. The signal of this system is approximately $3$, hence the lower noise level. 

Each prior knowledge setting had an associated model search space for which AIM was applied. The data was separated into environments, in each of which all but the inhibited species evolved over time. AIM was applied to each environment individually, and the resulting subnetworks were averaged to produce the full network estimates.

\subsubsection{Results}
\begin{figure}[h]
  \centering
    \includegraphics[width=\textwidth]{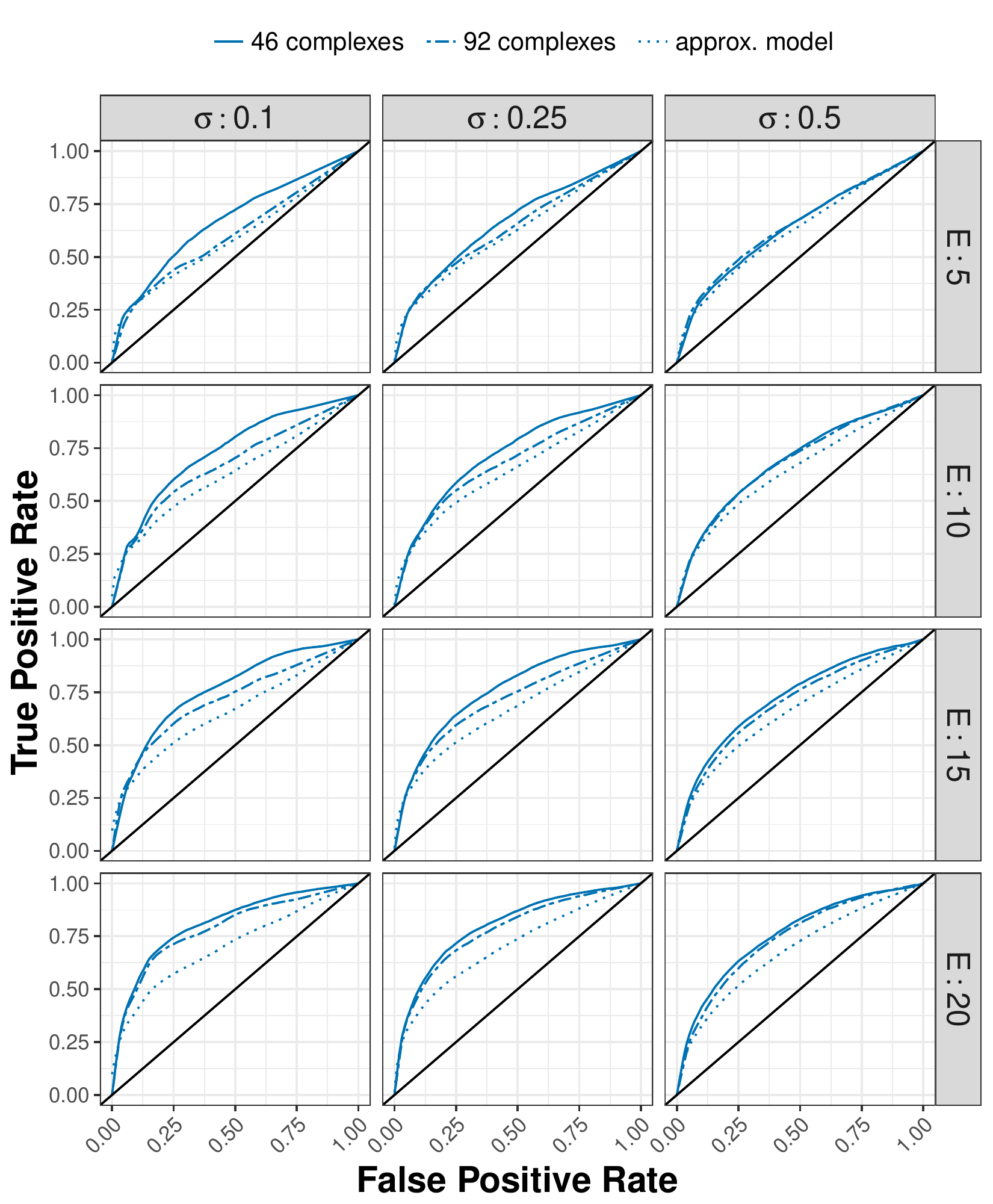}
  \caption{Pointwise average of the ROC curves, stratified according to noise level and number of environments.}  
  \label{fig:simC_rocs}
\end{figure}

\begin{figure}
  \begin{tabular}{c}
    \includegraphics[width=.6\textwidth]{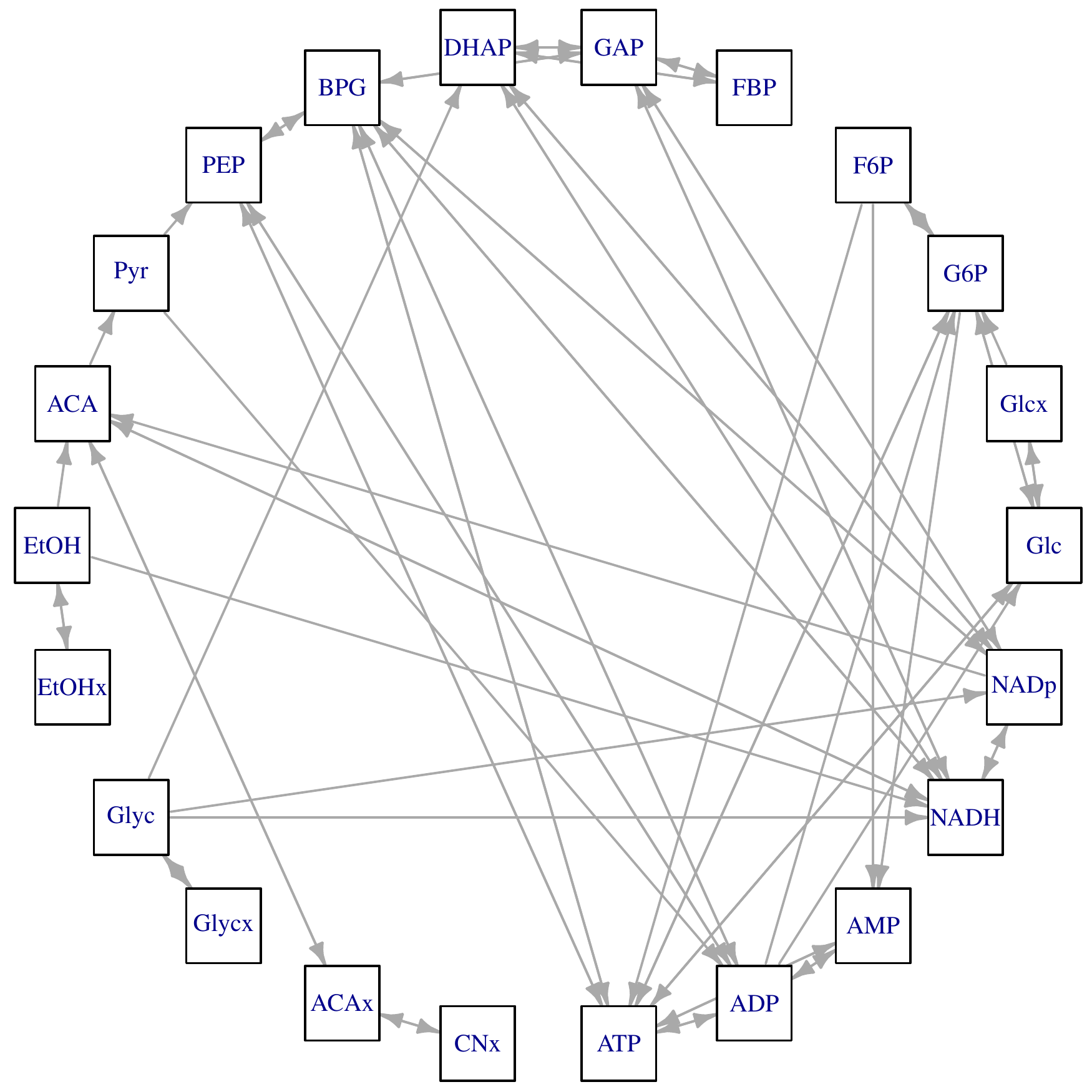} \\ 
    \includegraphics[width=.6\textwidth]{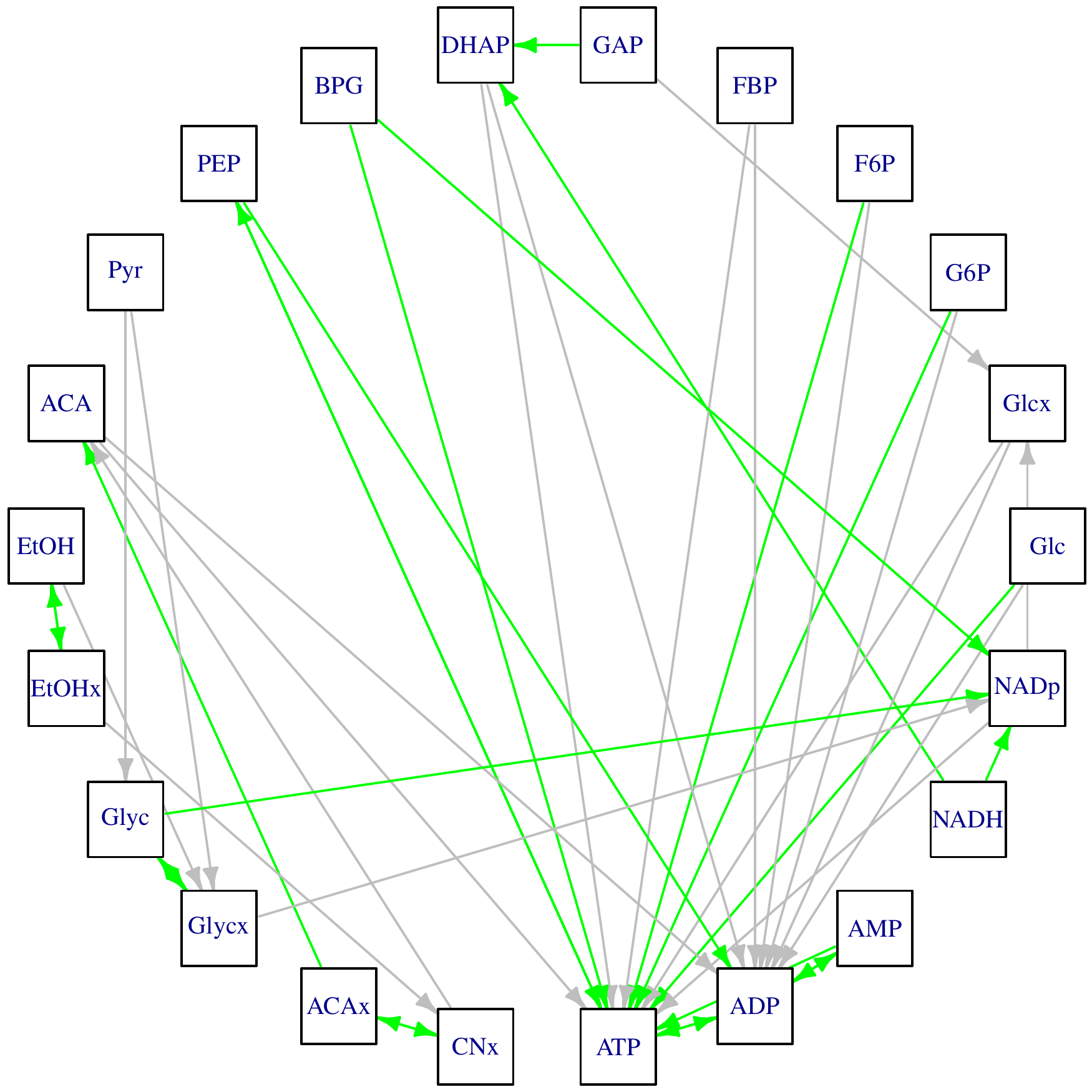} 
  \end{tabular}
\caption{The true glycolysis network (upper). The graph with the two most reported children for each node (lower). True edges are green. $E = 20$, $\sigma = 0.25$ and second setting with superset of complexes known.}
\label{fig:netw_summary}
\end{figure}

The ROC curves for the network estimator were calculated for each of the 100 replications. The average curves are in Figure \ref{fig:simC_rocs}. Not surprisingly the performance decreased with increasing noise, but more importantly we see a clear improvement with the number of environments. The estimated network and the true network are summarised in Figure \ref{fig:netw_summary}. We note that the approximate model has the overall worst performance in terms of network recovery, while the two models that incorporate prior knowledge by restricting the search space perform better. Though we do identify aspects of the network reasonably well, it is also evident that there is room for improvement, especially when no prior knowledge is used.

\clearpage

\section{Simulation Studies}
\label{sec:MAK}
In this section we return to the mass action kinetics systems:
\begin{equation}
 \label{eq:mak_matrix2}
 \frac{dx}{dt} = (B-A)^T\mathrm{diag}(x^A)k, \qquad x(0) = x_0,
\end{equation}
where $k=(k_r)_{r=1}^R$ and $x^A=(x^{A_r})_{r=1}^R$, and $R\in\mathbb{N}$ denotes the number of reactions. When the stoichiometric matrices $A$ and $B$ are either not known at all or only partially known, we seek to identify them from a larger set of candidate reactions. We test the performance of AIM in such a challenge through two simulation studies. 

\subsection{Simulation Study I}
In this section we compare AIM to another ODE network recovery algorithm GRADE, provided by \cite{Chen:2016}. GRADE is a nonparametric inverse collocation method. It replaces a parametric form of $f$ with a basis function expansion assuming an additive form, i.e., any given coordinate of $f$ depends on the other coordinates in an additive manner. GRADE was shown quite effective in simulation studies and applications by  \cite{Chen:2016}.

\subsubsection{Simulation study design}
The setting of this simulation study is a replicate of the simulation study in Section 5.3 in \cite{Chen:2016}. We consider five independent Lotka-Volterra systems, i.e., for $k=1,...,5$ we let
\begin{equation}
 \label{eq:lotka-volterra}
 \begin{aligned}
  \frac{dx_{2k-1}}{dt} &= 2x_{2k-1}(t) - v x_{2k-1}(t) x_{2k}(t) \\
  \frac{dx_{2k}}{dt} &= v x_{2k-1}(t) x_{2k}(t) - 2x_{2k-1}(t).
 \end{aligned}
\end{equation}
Note that the above ODE can be cast as a mass action kinetics system with $10$ species and $15$ reactions.

For each of the $E$ environments we drew the initialisation uniformly at random from $[0,4]$ and solved \eqref{eq:lotka-volterra} for $t\in [0, 5]$. Observations were extracted at $t=0,0.1,0.2,...,5$ with additive Gaussian noise. AIM was applied with a single linear interpolation smoother and GRADE used a spline smoother for smoothing the data and a monomial basis expansion of size 3 in \eqref{eq:approx_loss_integral_np}.

AIM searched ODE solutions using mass action kinetics reactions on the form
\begin{equation}
 \label{eq:mak_grade}
 \begin{array}{rcll}
  X_i + X_j &\rightarrow& 2X_i, &\qquad i,j = 1,...,10, \ i\neq j. \\
  X_i &\rightarrow & 2X_i,& \qquad i = 1,...,10. \\
  X_i &\rightarrow & 0,& \qquad i = 1,...,10. 
 \end{array}
\end{equation}
The search space thus consisted of $p=110$ reactions.

The following simulation parameters were used:
\begin{center}
\begin{tabular}{ccc}
\begin{tabular}{c}
Parameter \\
\hline
$E$ \\
$v$ \\
$\sigma$ 
\end{tabular}
& 
\begin{tabular}{cccc}
\multicolumn{4}{c}{Values} \\
\hline
2 & 4 & 8 & \\
1 & 3 & 5 & 7 \\
0.5 & 1 & 2 &
\end{tabular}
& 
\begin{tabular}{l}
Description \\
\hline
Number of environments \\
Interaction parameter \\
Standard deviation of additive noise 
\end{tabular}
\end{tabular}
\end{center}
The noise level was intentionally relatively large, as this ODE system is far easier to recover than those of the other systems considered in this paper. Each simulation was replicated 100 times.

\subsubsection{Results}
The ROC curves were derived for each simulation setting and method. A summary of the ROC curves is presented in Figure \ref{fig:grade_rocs_single} for $v=5$. Similar summaries for the remaining values of $v$ are found in the supplementary material. 

\begin{figure}[h]
  \centering
    \includegraphics[width=\textwidth]{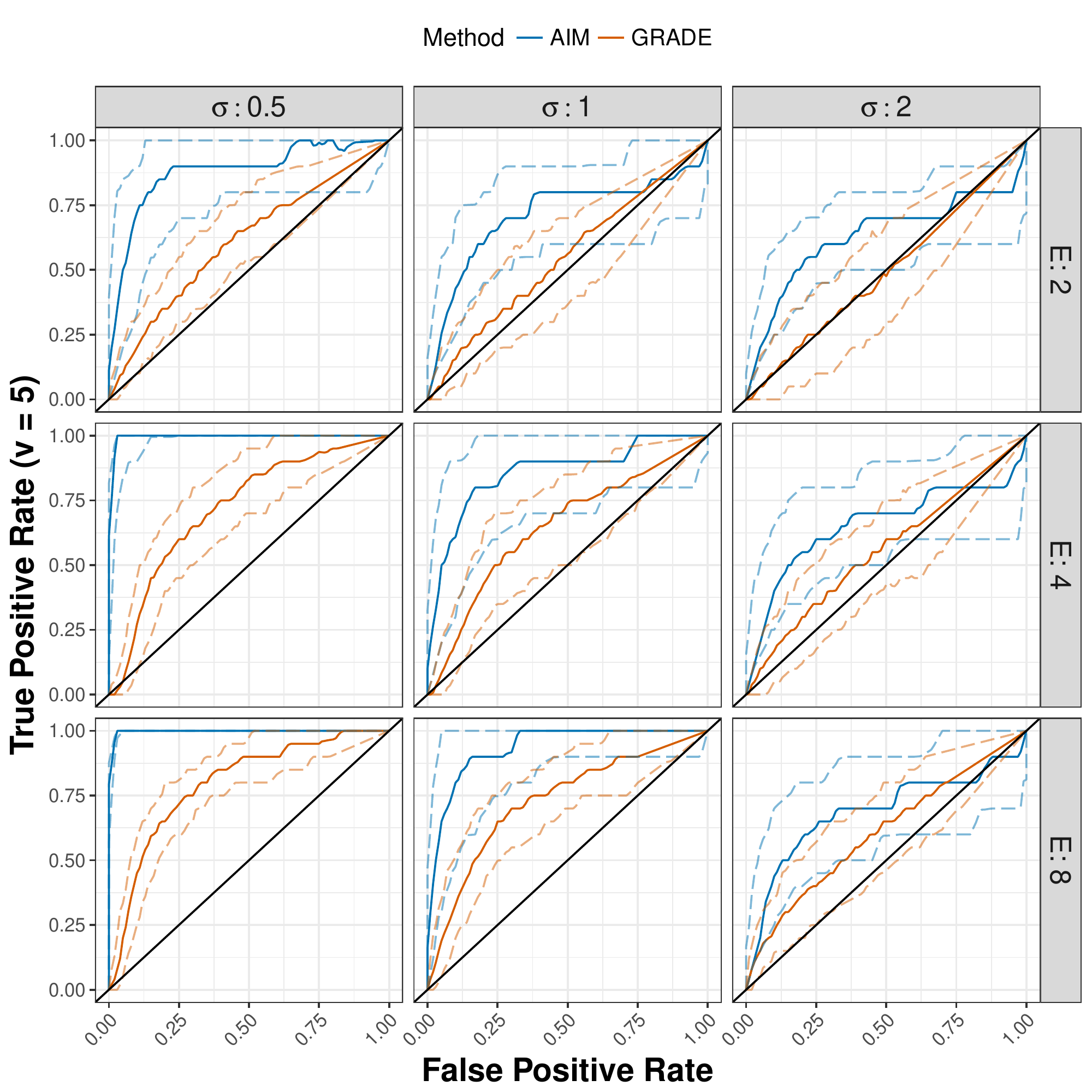}
  \caption{Pointwise median and $5^{\mathrm{th}}$ and $95^{\mathrm{th}}$ percentiles of ROC curves for the Lotka-Volterra system with $v=5$, stratified according to noise level and number of environments.}  
  \label{fig:grade_rocs_single}
\end{figure}

Across all noise levels, number of environments and interaction parameters we see that AIM generally performs better than GRADE. We ascribe this to the additivity assumption in GRADE, as we see improvements for decreasing $v$. Surprisingly, AIM works to an acceptable degree even for $\sigma=2$, which corresponds to a signal-to-noise ratio of $1$.

\clearpage

\subsection{Simulation Study II}
\label{sec:simB}
In this section we report the results from an extensive simulation study, whose purpose was to quantify how well AIM (in its concrete form of Algorithm \ref{alg:AIM_wrap}) identifies the correct reaction network. 

\subsubsection{Estimators} \label{sec:simest}
AIM was compared to an exhaustive gradient matching method (EGM), inspired by \cite{Babtie:2014}. See Appendix \ref{sec:comp_asp} for details on its implementation. Even though it relies on an inverse collocation method, this approach is computationally expensive as it selects the reactions based on best subset selection applied to each species separately. This computer intensive inverse collocation method attempts to get the most information out of the approximate loss function, by finding global minima on lower dimensional subspaces. 

Solving the best subset selection problem for each species separately only induces the global best subset selection solution if each coordinate of $\theta$ induces a single edge in the network. This property holds for linear ODE systems and does not hold for most mass action kinetics systems. This simulation study restricts the attention to reactions on the form $X_i + X_j \rightarrow 2X_i$, $i \neq j$, hence each reaction corresponds to a bidirectional edge between node $i$ and $j$ -- as well as a self-edge in both nodes. Thus, in this particular simulation study, EGM will provide the same networks as a best subset selection performed over all possible reactions. EGM was not applied to the examples considered in Section \ref{sec:app}, as the number of species was too large for an exhaustive search to be computationally feasible. 

Each method reported estimated reaction networks consisting of up to $5d$ reactions. EGM used the Gaussian process smoother described in \cite{Babtie:2014}, IM used a linear interpolation smoother and AIM used both smoothers. In order to produce additional initialisations for AIM, the integral matching estimates from each smoother were produced with and without standardising the coordinates of the process. 
Both AIM and IM used the elastic net penalty (\cite{Zou:2005}) with $\alpha = 0.25$. 


\subsubsection{Simulation study design} \label{sec:simdesign}
Data was drawn from reaction networks composed of reactions on the form:
\begin{equation}
 \label{eq:mak_enzyme}
 X_i + X_j \rightarrow 2X_i
\end{equation}
where $i,j = 1,...,d$ and $i\neq j$, with a total number of reactions at $p=d(d-1)$. Time course data from $E$ environments were drawn. Each environment was given its own initial condition produced as follows: between one and four distinct chemical species were selected at random. In each environment all but the selected species were knocked down by $50\%$ from their equilibrium value and the initial abundance of the selected species were increased by the total mass knocked down. The initial conditions were rescaled to have an average of 5. Since the total number of molecules is preserved by reactions on form \eqref{eq:mak_enzyme}, the average signal strength is approximately 5. 

Data were sampled at $t=0$ and $t=2^{i/2}$, for $i=-5,-4,...,2,3$, all with additive Gaussian noise. Each species $i=1,...,d$ was given $\alpha=1,2$ true reactions, i.e., a total of $d\alpha$ reactions on the form \eqref{eq:mak_enzyme} had rate parameter $1$ and the remaining $0$.

The following simulation parameters were used:
\begin{center}
\begin{tabular}{ccc}
\begin{tabular}{c}
Parameter \\
\hline
$d$ \\
$\alpha$  \\
$E$ \\
$n$ \\
$\sigma$ 
\end{tabular}
& 
\begin{tabular}{ccc}
\multicolumn{3}{c}{Values} \\
\hline
7 & 9 & 11 \\
1 & 2 \\
2 & 4 & 8  \\
10 \\
0.1 & 0.5 & 1 
\end{tabular}
& 
\begin{tabular}{l}
Description \\
\hline
Number of species \\
Number of true reactions per species  \\
Number of environments \\
Number of data points per environment \\
Standard deviation of additive noise 
\end{tabular}
\end{tabular}
\end{center}
For each combination of the simulation parameters, 100 replications of the above simulation experiments were conducted.

\subsubsection{Results}
We first report the recovery of the true network. For each replicate and simulation parameter combination the receiver operating characteristic (ROC) curves of the network were derived. Pointwise averages over the replicates are illustrated in Figure \ref{fig:roc_single} for $d=9$ and $\alpha=1$. The remaining curves can be found in the supplementary material.

\begin{figure}[h]
  \centering
    \includegraphics[width=\textwidth]{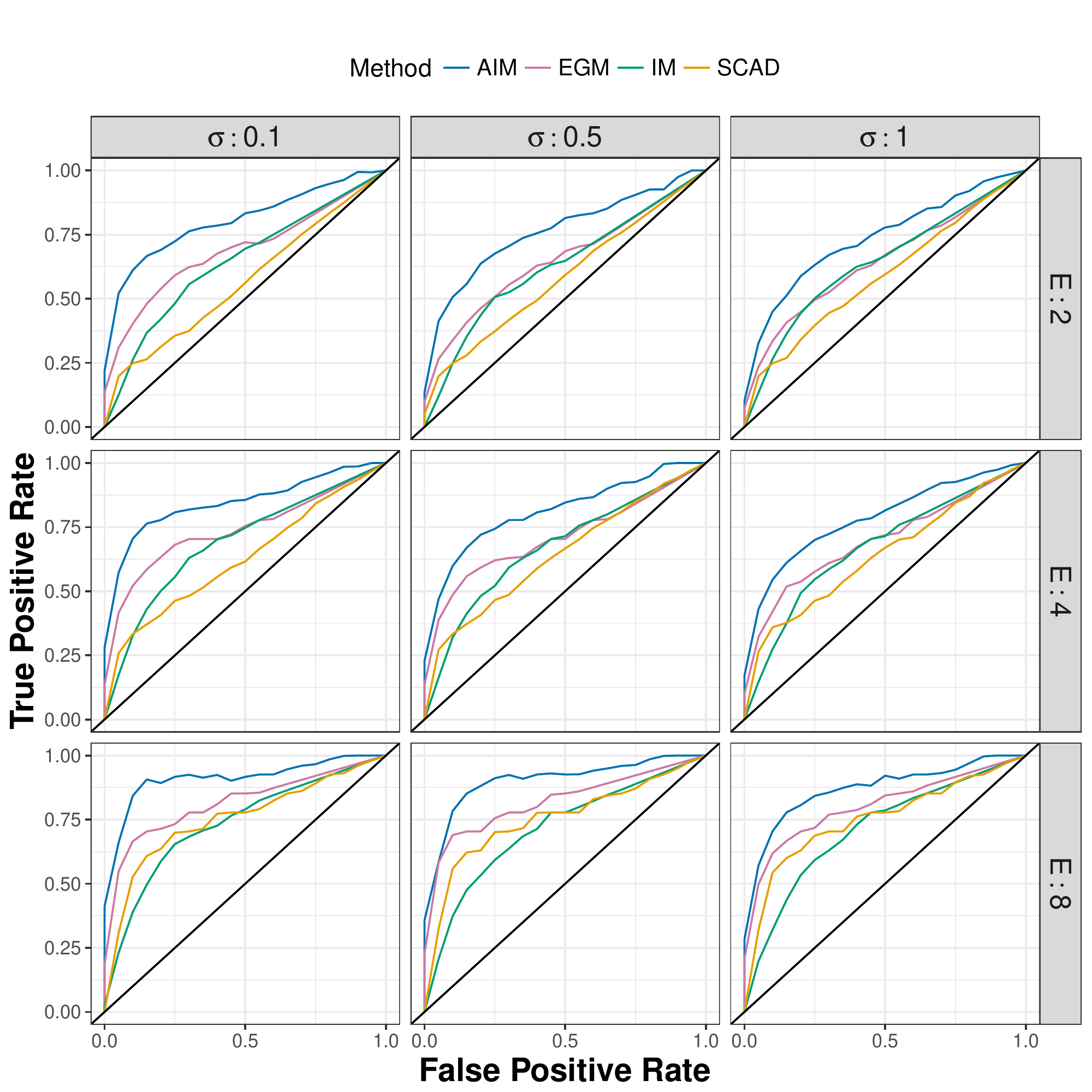}
  \caption{Pointwise averaged ROC curves of the network estimates for $d=9$ and $\alpha=1$, stratified according to noise level and number of environments.}  
  \label{fig:roc_single}
\end{figure}

\begin{figure}[h]
  \centering
    \includegraphics[width=\textwidth]{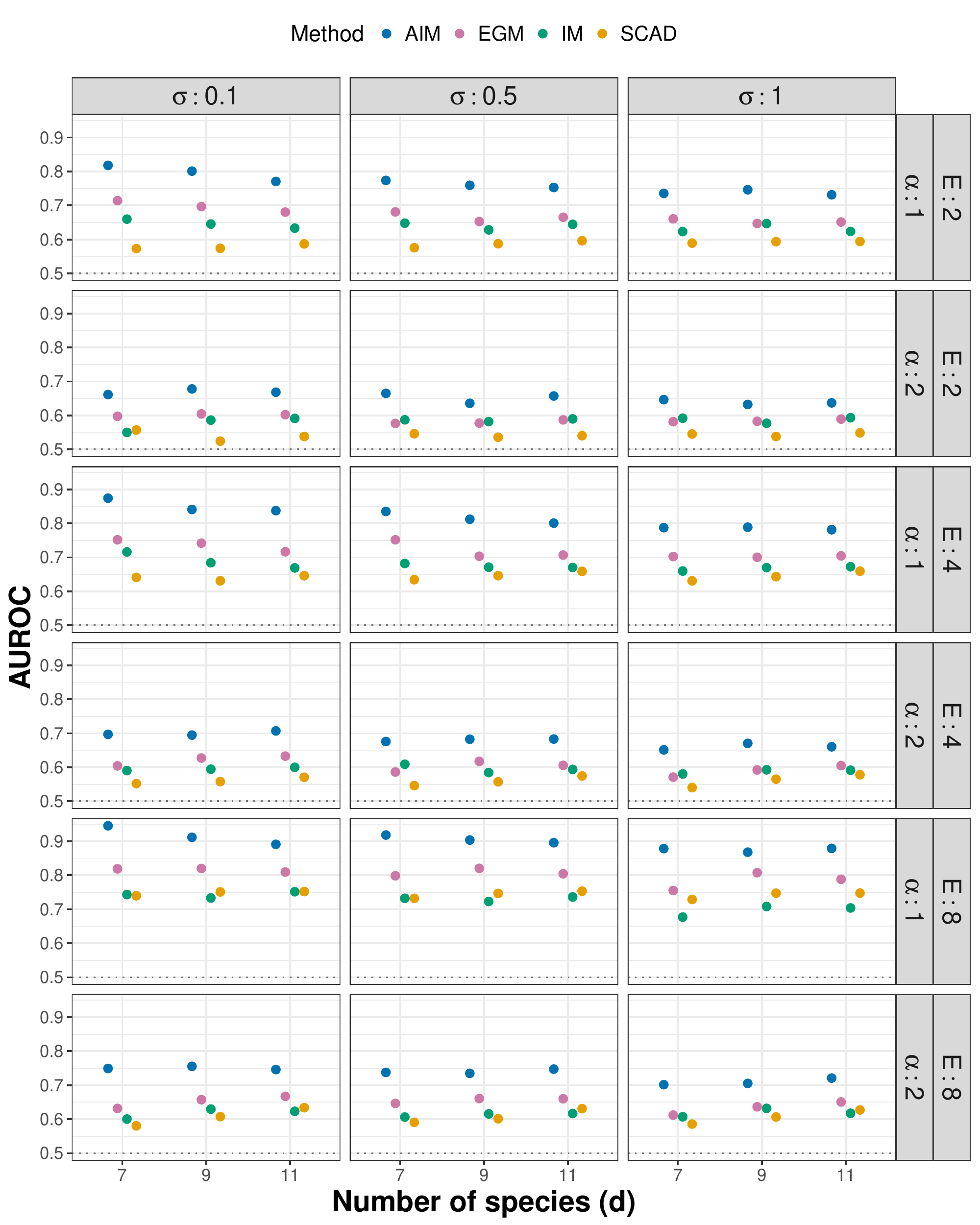}
  \caption{Median AUROC across the 100 replications and stratified according to the simulation settings.}  
  \label{fig:auroc}
\end{figure}

From Figure \ref{fig:roc_single} we see that AIM consistently recovered the network better than the other methods. IM and SCAD were the worst performing methods with SCAD improving the most with increasing number of environments, though not reaching the level of EGM and AIM. 

These tendencies are repeated in the other figures in the supplementary material, with an overall decrease in performance for $\alpha=2$. Figure \ref{fig:auroc} provides an overview of the area under the ROC curves (AUROC) across simulation settings.  AIM generally had the largest median AUROC values across all settings. For all methods we also see improvements when increasing the number of environments and that increasing the number of species for most scenarios decreases the performance. Generally, for all methods the network recovery performances dropped considerably for $\alpha=2$.

\begin{figure}[h]
  \centering
    \includegraphics[width=\textwidth]{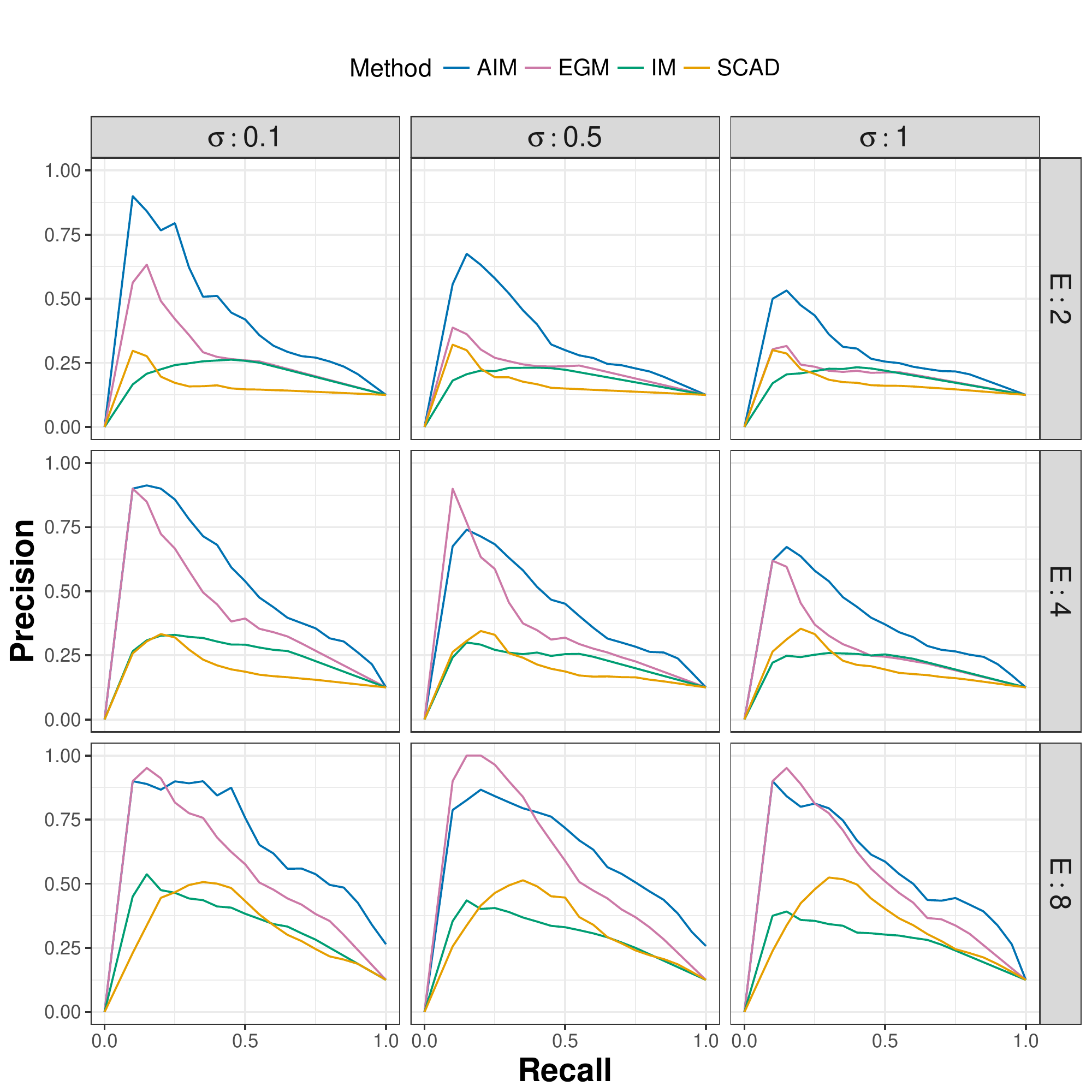}
  \caption{Pointwise averaged precision-recall curves of the reactions recovered for $d=9$ and $\alpha=1$, stratified according to noise level and number of environments.}  
  \label{fig:precis_recall_single}
\end{figure}

Next we report the recovery of the true reactions. We visualise their performance by their precision-recall curves. In Figure \ref{fig:precis_recall_single} the pointwise averaged precision-recall curves for $d=9$ and $\alpha=1$ are presented. The remaining curves can be found in the supplementary material.

From Figure \ref{fig:precis_recall_single} we see that EGM recovered most correct reactions early in the recovery for $E$ large. But after recovering 20--35\% of the true reactions AIM surpassed EGM in reaction recovery performance. All methods improved considerably with increasing number of environments. These results match what we observed for the network recovery to some degree.

The network ROC curves and the reaction precision-recall curves together suggest that EGM recovers the first few reactions and network edges accurately, but AIM is more accurate when more reactions are reported.
\clearpage

\begin{figure}[h]
  \centering
    \includegraphics[width=\textwidth]{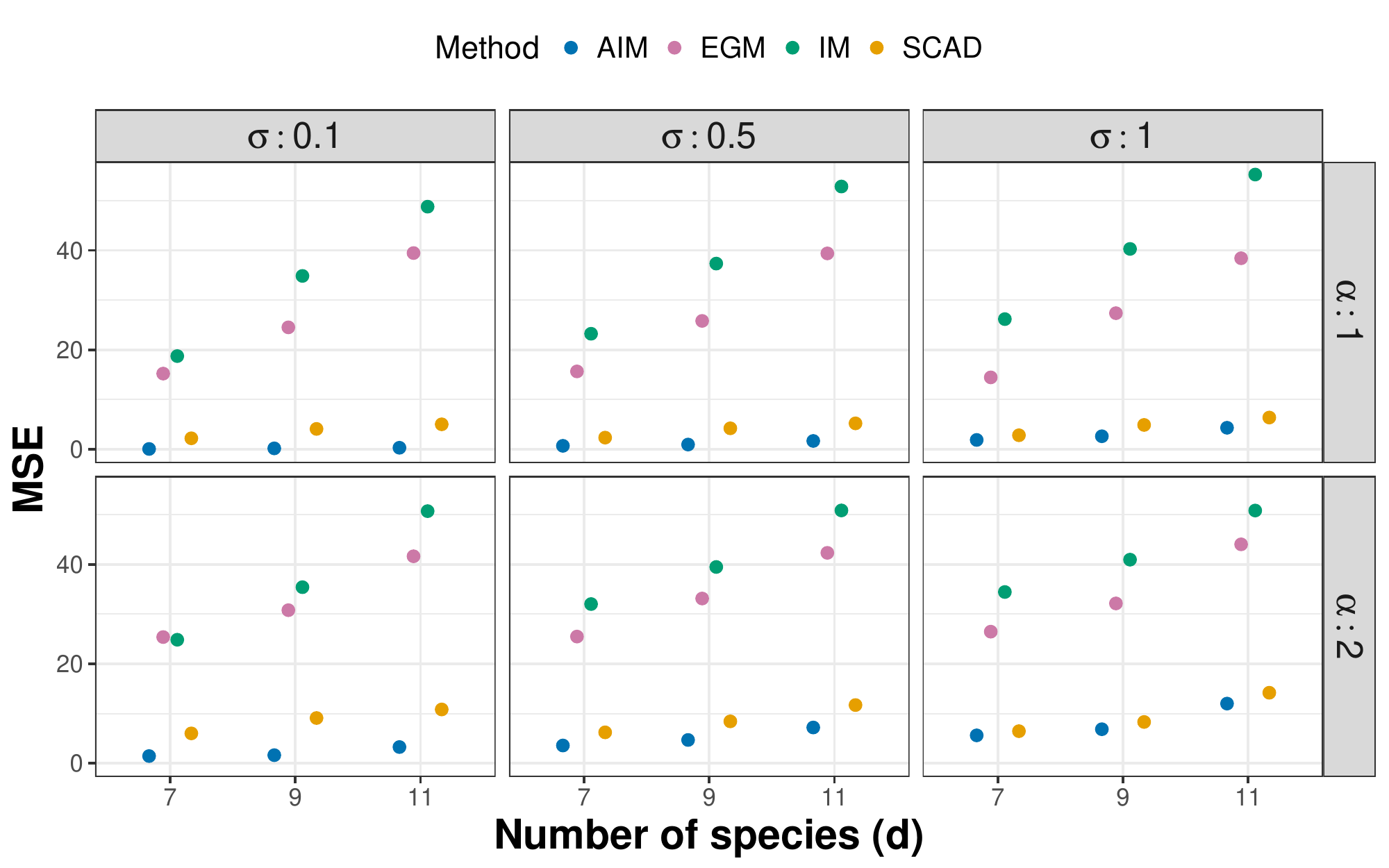}
  \caption{Medians of the mean squared error between the tuned trajectories and the true trajectory for $E=4$.}  
  \label{fig:mspes}
\end{figure}

The mean squared error of the estimated reaction networks were also assessed. A single model was selected for each method by minimising the mean squared error on an independent test set. The squared error between the trajectories produced by the selected model and the true trajectory at each time point was derived. Medians of the mean squared error are presented in Figure \ref{fig:mspes}. 

We see that the methods using the ODE-based loss (AIM and SCAD) have much smaller mean squared error than the methods based on the approximate loss. That the mean squared error is so large for IM and EGM can be explained as a bias phenomenon similar to the one observed for the Michaelis-Menten example as illustrated in Figure \ref{fig:MM}. IM without a penalty and a linear interpolation smoother -- as used in this simulation study -- is expected to be relatively unbiased but with a large variance. However, the sparsity enforcing penalty introduced an additional bias, and the resulting trajectories of the fitted ODE did not match the truth very well in general (data not shown). For EGM the conclusion is the same, but the argument is the other way around. This method used a Gaussian process smoother, which should decrease the variance of the parameter estimates, but the under-estimated slopes introduced a stronger bias. Again, the resulting trajectories of the ODE fitted using EGM did not match the truth very well. 
 Though the mean squared error suggests that the approximate loss functions provide quantitatively incorrect estimates, we did find qualitatively correct network and reaction recovery for those methods.

\begin{figure}[h]
  \centering
    \includegraphics[width=0.9\textwidth]{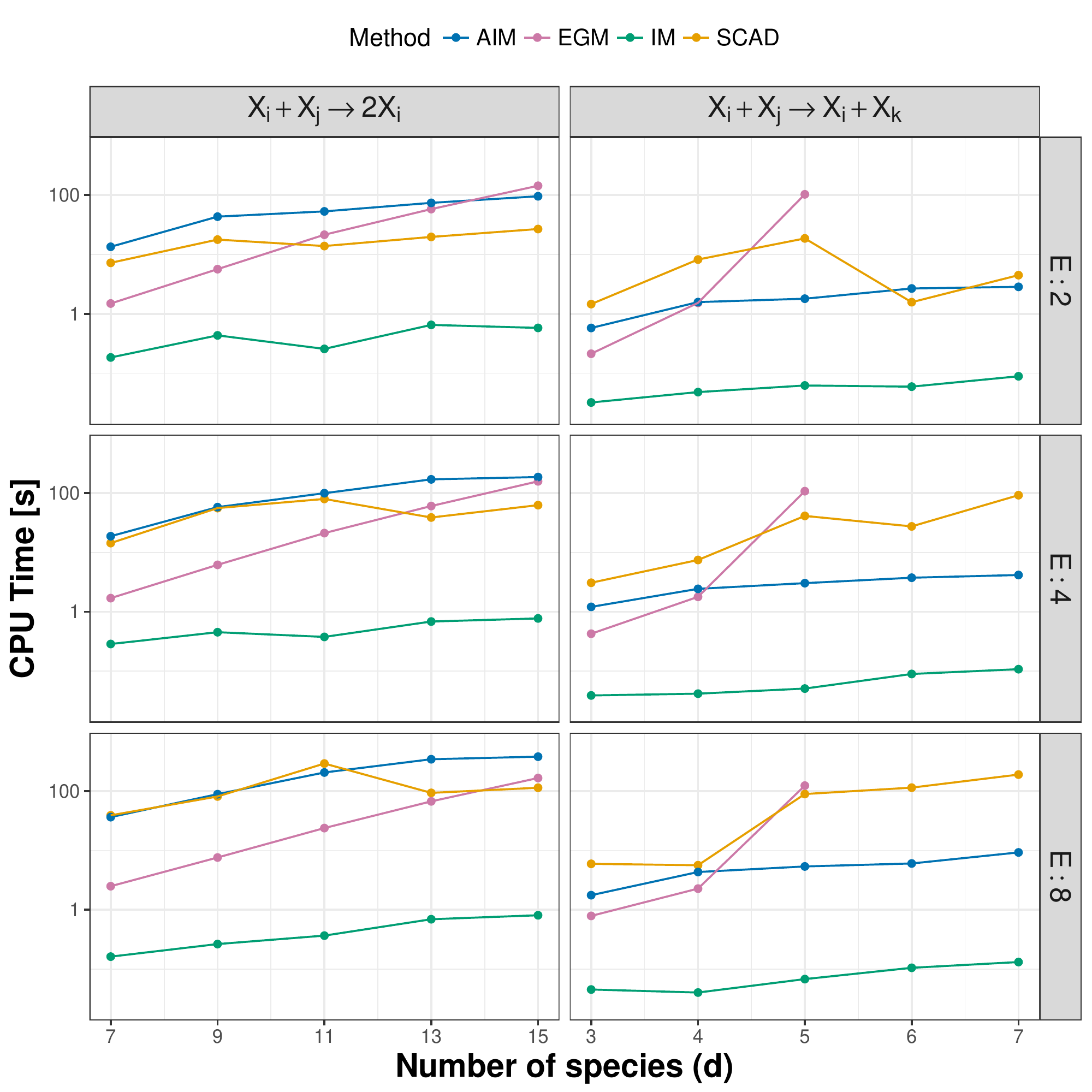}
  \caption{Median computation time for the two reaction collections. $\sigma = 0.5$ and $\alpha=1$.}  
  \label{fig:cpus}
\end{figure}

Finally we report computation times. The median computation time over 10 replications can be found in Figure \ref{fig:cpus} for two collections of reactions: $X_i + X_j \rightarrow 2X_i$, $i\neq j$ and $X_i + X_j \rightarrow X_i + X_k$, $j\neq i\neq k$. The model search space size of the latter grows faster with the number of species and it quickly becomes a challenge for EGM. In fact EGM was excluded for $d>5$ due to infeasible computation times. For  $d$ small, AIM is somewhat slow, however in terms of scalability with $d$ AIM resembles IM more than EGM. 


%
%

\clearpage

\section{Discussion}

Collocation based estimation of parameters in ODE systems is computationally less demanding than the least squares method that relies on numerical solutions of the ODE systems. \cite{Calderhead:2009} demonstrated this in a Bayesian setting and proposed a collocation method based on Gaussian processes, \cite{Babtie:2014} relied on a Gaussian process collocation based search of the model space similar to the EGM method that we have implemented, and \cite{Chen:2016} relied on penalised  collocation based estimation for their method GRADE (Graph Reconstruction via Additive Differential Equations). Based on these and similar results we sought to develop a scalable inference framework for mass action systems, but we found several challenges, and the present paper represents a synthesis of how we dealt with these challenges. We discuss below how the most important challenges were addressed.

\subsection{Bias} We found that penalised collocation methods were computationally fast, but even if they did recover qualitatively the correct network and reactions for realistic data sizes, the resulting parameter estimates were biased.  The bias was induced partly by the initial smoothing and partly by the penalisation, and the fitted model would not reproduce very accurately the solution trajectories of the true data generating ODE system. Moreover, the results would be sensitive to the precise choice of initial smoother. We found that among the collocation methods, our proposed integral matching (IM) estimator obtained by minimising \eqref{eq:approx_loss_aim} has reasonable statistical properties. 

\subsection{Penalised least squares} To test if penalised least squares methods are feasible for large systems we implemented a number of algorithms for numerical minimisation of the penalised least squares loss including the proximal gradient algorithm with screening as presented in Appendix \ref{sec:screen}. Sparsity and screening combined with fast solvers of the sensitivity equations makes it possible to apply these algorithms even for fairly large systems. However, the sparsity inducing penalty still induces a bias of the resulting estimates, which can also be quite dependent on the initialisation of the optimisation algorithm due to local minima of the objective function. We illustrated that least squares with the SCAD penalty achieved rather accurate estimates in terms of mean squared error from the true trajectories, but in terms of network recovery it was inferior to the other methods considered -- in particular IM, which is much faster.  

\subsection{Parameter scale} A parameter in an ODE system typically controls the rate of a reaction, and the bias induced by the penalty results in reaction rates being underestimated. It is our experience that the bias induced by the penalty can have quite substantial effects for nonlinear ODE systems, and the choice of parameter scale determines this bias together with the combined effect of the penalty term. Standardisation as used in regression models, e.g. in the R package glmnet \citep{Friedman:2010}, for bringing the parameters on a common scale is not directly applicable. We suggest adaptive rescaling as given by \eqref{eq:adapt_sc}, which does require a pilot estimate of the unknown parameter unless $f$ is linear in $\theta$. However, we did not find this to be a data-driven panacea for the choice of parameter scale, and we ended up concluding that the unpenalised estimator given by \eqref{eq:loss_no_pen} had better statistical properties in our experiments. 

\subsection{Combining methods} Our combined AIM algorithm uses the fast collocation method IM to obtain good initialisation parameters for the least squares method. Moreover, AIM in the form of Algorithm \ref{alg:AIM_wrap}  -- which we have extensively tested -- uses multiple smoothers to achieve an even greater variety of initialisations, and it introduces sparsity in the least squares method by restricting the parameter space. The stratified ranking was proposed as a way to aggregate the resulting models into a sequence of models indexed by the number of nonzero parameters. Clearly, alternative aggregations are possible, e.g. using a weighted average.
Moreover, the simulation study in Section \ref{sec:simB} found that EGM performed slightly better than AIM for the first couple of reactions. As EGM performs the first couple of search iterations fairly quickly, a hybrid approach for initialisation suggests itself using EGM for the first couple of reactions and IM for the remaining reactions. We have not investigated if alternative aggregation schemes or hybrid approaches for initialisation could further improve the statistical properties of the algorithm.

\subsection{Network recovery} We demonstrated that AIM has good network recovery properties in a number of different examples and compared to several alternatives. In the Lotka-Volterra example it was, for instance, demonstrated that AIM was far superior to GRADE \citep{Chen:2016}. This is perhaps unsurprising given that GRADE assumes an additive form of $f$, but we emphasise this to argue that additivity is a quite strong assumption, which is unlikely to hold for many ODE systems of practical relevance. 

AIM also performed well in the recovering of the \emph{in silico} network of protein phosphorylation, and it was superior in terms of AUROC to IM and SCAD considered in this paper as well as the two ODE-based solutions that participated in the eighth DREAM challenge. We did not participate in the challenge, but AIM would have been ranked second among all participants. We note that the top-ranked submissions including the winning team did not rely on a mechanistic model -- the submission only required network edge weights. The winning team constructed edge weights via tests for nonlinear functional relations without the constraints of an ODE system, and were in this way better able to capture the correct network structure (see Supplementary material on Team 7 in \cite{DREAM8}). However, such methods are not capable of predicting e.g. intervention or perturbation effects. It is clearly of interest to utilise such network estimates as prior information for learning ODE systems, and we demonstrated how this can be done in our framework for the discovery of the glycolysis network. For the DREAM challenge it would make an unfair comparison if we were to piggyback on the published top-ranked network for this particular data set, hence we ran AIM in this example without any prior network restrictions.  

\subsection{Conclusion} The AIM algorithm was presented and demonstrated to have good statistical properties for realistic data structures and sizes. The implementation of AIM also demonstrated that it is possible to learn large ODE systems via least squares methods -- even if this is computationally heavy. Further improvements may be possible, e.g. to account for a more complicated noise structure than additive, uncorrelated noise. In the light of the linear noise approximation, described in detail by \cite{Wallace:2011}, the noise in mass action kinetics systems scales with the signal. We have partially addressed this by rescaling the observation weights as given by \eqref{eq:adapt_weight}, which will adjust the weights according to the variance of each species. However, we have not investigated ways to adjust for a more complicated variance structure.  

Our intention with AIM and the associated R package \emph{episode} is to provide a thoroughly tested, applicable and useful framework for learning ODE systems using state-of-the-art methods. This should be of use to experimentalists, and it should be able to serve as a benchmark for further developments. The R package currently supports polynomial and rational systems in certain parameterisations, and it is implemented in a modular way that allows for easy addition of new parameterised families of ODE systems. Doing so, the entire framework consisting of data structures, ODE solvers and optimisers including AIM are then directly available.

\clearpage
\appendix

\section{Computational Aspects}
\label{sec:comp_asp}

\subsection{Sensitivity equations and approximative gradients}
Let $x:\mathbb{R}\rightarrow \mathbb{R}^d$ solve the ODE
\begin{equation}
 \label{eq:ode_}
 \frac{dx}{dt} = f(x, \theta), \qquad x(0) = x_0.
\end{equation}
The derivative of $x$ with respect to $\theta\in\mathbb{R}^p$, i.e., the matrix valued function $x_{\theta}:\mathbb{R}\rightarrow \mathbb{R}^{d\times p}$, solves another ODE system:
\begin{equation}
 \label{eq:ode_sens_theta}
 \frac{dx_{\theta}}{dt} = \frac{\partial f}{\partial x}(x, \theta)x_{\theta} + \frac{\partial f}{\partial \theta}(x, \theta), \qquad x_{\theta}(0) = \mathbb{O}_{d\times p},
\end{equation}
where $\mathbb{O}_{d\times p}$ is the $d\times p$-dimensional zero-matrix. Analogously, the derivative of $x$ with respect to $x_0$,  $x_{x_0}:\mathbb{R}\rightarrow \mathbb{R}^{d\times d}$, solves the ODE:
\begin{equation}
 \label{eq:ode_sens_x0}
 \frac{dx_{x_0}}{dt} = \frac{\partial f}{\partial x}(x, \theta)x_{x_0}, \qquad x_{x_0}(0) = \mathbb{I}_d,
\end{equation}
with $\mathbb{I}_{d}$ the $n$-dimensional identity matrix. The equations \eqref{eq:ode_sens_theta} and \eqref{eq:ode_sens_x0} are called the \textit{sensitivity equations} of \eqref{eq:ode_}. Notice that once the original system is solved, the columns of the sensitivity equations can be solved independently.

The sensitivity equations are often solved simultaneously with the original system \eqref{eq:ode_}. Even if \eqref{eq:ode_} requires a computationally intensive solver (e.g., a solver with adaptive step length or implicit solvers), the sensitivity equations often only require simple solvers like the Euler scheme to be accurate. There are two reasons for this. Firstly, the sensitivity equations are (time-inhomogeneous) affine ODE systems which are often less sensitive. Secondly, the exact gradient is not necessary to optimise a smooth function -- an approximate gradient pointing in roughly the same direction will suffice.

A method for deriving even faster approximate gradients to 
\begin{equation}
 \label{eq:loss_eysm}
 \ell_y \coloneqq \sum_{t \in \mathcal{C}}{\|y_t - \int_0^t{f(x(s, \theta), \theta) ds}\|_2^2}
\end{equation}
was proposed by \cite{Mikkelsen:2015} and inspired by inverse collocation methods. It goes as follows: assuming that $\theta_0$ is the current value of $\theta$ in the optimisation, then minimise
\begin{equation}
 \label{eq:loss_eysm_approx}
 \theta \mapsto \sum_{t \in \mathcal{C}}{\|y_t - \int_0^t{f(x(s, \theta_0), \theta) ds}\|_2^2}
\end{equation}
to produce the next step. Though these approximate solutions are not guaranteed to improve the original loss function, they still produce fast and approximate descent directions. If the approximate solution does not improve the loss function, it is suggested to take one classic gradient-based step before retrying the approximate solution. 

The above approach is equivalent to using the Gauss-Newton method on the original loss function, but ignoring the first term of the right hand side of \eqref{eq:ode_sens_theta}, when calculating the differentials.

\subsection{Proximal gradient and screening methods}
\label{sec:screen}
The penalised ODE loss function
\begin{equation}
 \label{eq:loss_}
 \ell_y(\theta) \coloneqq \frac{1}{2}\sum_{i=1}^{n}{ \sum_{l=1}^{d} { w_{i,l}(y_l(t_i) - x_l(t_i,\theta))^2 } } +\lambda\sum_{j=1}^{p}{v_j \mathrm{pen}(\theta_j)},
\end{equation}
is optimised using a proximal gradient method, as described in \cite{Hale:2007} for $\ell^1$ penalties. For non-convex penalties, like SCAD and MCP, this method is combined with the majorisation method by \cite{Fan:2001}. 

The proximal gradient method for \eqref{eq:loss_} thus starts with initialisation $\theta^0$ and the proceeds with
\begin{equation}
 \label{eq:prox_alg}
 \theta^{k+1} = \mathrm{prox}(\theta^k, p^k, \tau), \qquad \text{for } k = 0,1,...
\end{equation}
where the proximal operator is defined as
\begin{equation}
 \label{eq:prox}
 \mathrm{prox}(\theta, p, \tau) \coloneqq \mathrm{sign}(\theta - \tau p) \circ \mathrm{max}(0, |\theta - \tau p| - \lambda \mu(\theta)).
\end{equation}
The vector $p^k$ is the derivative of $\frac{1}{2}\sum_{i=1}^{n}{ \sum_{l=1}^{d} { w_{i,l}(y_l(t_i) - x_l(t_i,\theta))^2 } }$ at $\theta^k$, i.e.,
\begin{equation}
 \label{eq:diff_prox}
 p^k \coloneqq - \sum_{i=1}^{n}{ \sum_{l=1}^{d} { w_{i,l}(y_l(t_i) - x_l(t_i,\theta^k)) \frac{dx_l}{d\theta}|_{\theta = \theta^k}(t_i) } }
\end{equation}
and the sensitivity equation is solved using the approximative methods described above. For non-convex penalties, the majorisation amounts to replacing the penalty weights in each step by $v_j \circ \frac{d^2\mathrm{pen}}{d\theta_j^2}(\theta_j)$. 

In \eqref{eq:prox_alg} the step length $\tau$ is chosen through backtracking. Moreover, not all coordinates of $\theta$ changes in each step of \eqref{eq:prox_alg}. This is due to the sparsity inducing property of the proximal operator. In practice this means that many coordinates of the derivatives $p^k$ are calculated (using computationally intensive numerical solvers) and then never used. Computations are thus saved if occasionally the ODE system is screened for strong variables as follows: at every $n^{\mathrm{th}}$ step all coordinates of $p^k$ are evaluated. If $\theta^k = \mathrm{prox}(\theta^k, p^k, 1)$ (up to some numerical accuracy) then stop, else identify the active set $\mathcal{A}=\{i\mid \theta^k_i\neq 0 \text{ or } \theta^k_i \neq  \mathrm{prox}(\theta_i^k, p^k_i, 1) \}$ and run proximal gradient algorithm on $\mathcal{A}$ only until next screening.

\subsection{Exhaustive Gradient Matching}
Inspired by \cite{Babtie:2014} exhaustive gradient matching applies a best subset selection of parameters for explaining the dynamics of each chemical species individually. The individual results are then combined into a parameter estimate of the full ODE system. 

For an ODE system given by the field $f(x, \theta)$, then each coordinate of the solution satisfies 
\begin{equation}
  \label{eq:ode_single}
  \frac{dx_l}{dt} = f_l(x, \theta), \quad l=1,...,d
\end{equation}
where $f_l$ is the $l^{\mathrm{th}}$ coordinate of $f$ and $\theta\in \mathbb{R}^p$. Given smoothed curves $\hat{x}=(\hat{x}_l)_{l=1}^d$ for each coordinate, then the approximate inverse collocation loss function is
\begin{equation}
  \label{eq:collocation_loss}
  \ell(\theta)\coloneqq \frac{1}{2}\sum_{l=1}^d {\sum_{t \in \mathcal{C}}{\left(\frac{d\hat{x}_l}{dt}(t) - f_l(\hat{x}(t), \theta)\right)^2}}.
\end{equation}
If the field is linear in the parameters the above becomes the sum of squares,
\begin{equation}
  \label{eq:collocation_ss}
  \ell(\theta) = \sum_{l=1}^d { \|Y_l - X_l\theta \|_2^2 }.
\end{equation}
where $Y_l = \left(\frac{\hat{x}_l}{dt}(t)\right)_{t \in \mathcal{C}}$ and $X_l = \left(\frac{\partial f_l}{\partial \theta}(\hat{x}(t))\right)_{t \in \mathcal{C}}$ is a concatenation of the $\theta$-gradients of the field over the time points.

The exhaustive gradient matching method (EGM) goes as follows: for each $l=1,...,d$ construct $Y_l = \left(\frac{\hat{x}_l}{dt}(t)\right)_{t \in \mathcal{C}}$ and $X_. = \left(\frac{\partial f_l}{\partial \theta}(\hat{x}(t))\right)_{t \in \mathcal{C}}$. For any $\mathcal{K} \subseteq \{1,...,p\}$ let $X_l^{\mathcal{K}}$ denote the $\mathcal{K}$ columns of $X_l$ and let $\theta^{\mathcal{K}}\in\mathbb{R}^{|\mathcal{K}|}$. For $k=1,...,K$ find the subset $\mathcal{K}^l_k\subseteq \{1,...,p\}$ with $|\mathcal{K}^l_k|=k$ such that 
\begin{equation}
  \min_{\theta^{\mathcal{K}^i_{k}}} \frac{1}{2}\|Y_l - X_l^{\mathcal{K}_k}\theta^{\mathcal{K}_{k}}\|_2^2
\end{equation}
is minimal.

Each species now has a sequence of subsets of increasing size, $(\mathcal{K}^l_k)_{k=0}^K$. They are combined into a sequence of subsets representing the full system, $(\mathcal{K}_k)_{k=0}^{dK}$, by the union
\begin{equation}
 \mathcal{K}_k\coloneqq \bigcup_{i=1}^d\mathcal{K}_{\alpha_i(k)}^l.
\end{equation}
The $k$-dependent tuple $\alpha(k)=(\alpha_l(k))_{l=1}^d\in\{1,...,K\}^d$ is given by the recursion
\begin{equation}
 \alpha(0)=(0)_{l=1}^d, \qquad \alpha(k+1)=\alpha(k) + e_{l^*}, \quad k=0,...,Kd-1
\end{equation}
where the increments $e_{l^*}$ is $1$ at coordinate $l^*$ and zero elsewhere. The coordinate $l^*$ is chosen such that 
\begin{equation}
 \min_{\theta: j \notin \bigcup_{l=1}^d\mathcal{K}_{\gamma_i}^l \Rightarrow \theta_j=0} \ell(\theta), \qquad \gamma = \alpha(k) + e_{l^*}
\end{equation}
is minimal, i.e., the species whose next subset improves the loss the most determines the next full subset $\mathcal{K}_{k+1}$.

The EGM estimator becomes the best subset selection estimator of \eqref{eq:collocation_loss} if and only if each coordinate of the parameter vector $\theta$ affects only one edge in the network. This is the case for linear ODE systems and, if ignoring self-edges, also the case for the systems studied in Section \ref{sec:simB}. However, for most ODE systems a single coordinate often affects multiple network edges simultaneously.

\clearpage
\bibliography{bibliography}{}
\bibliographystyle{agsm}

\includepdf[pages={1-17}]{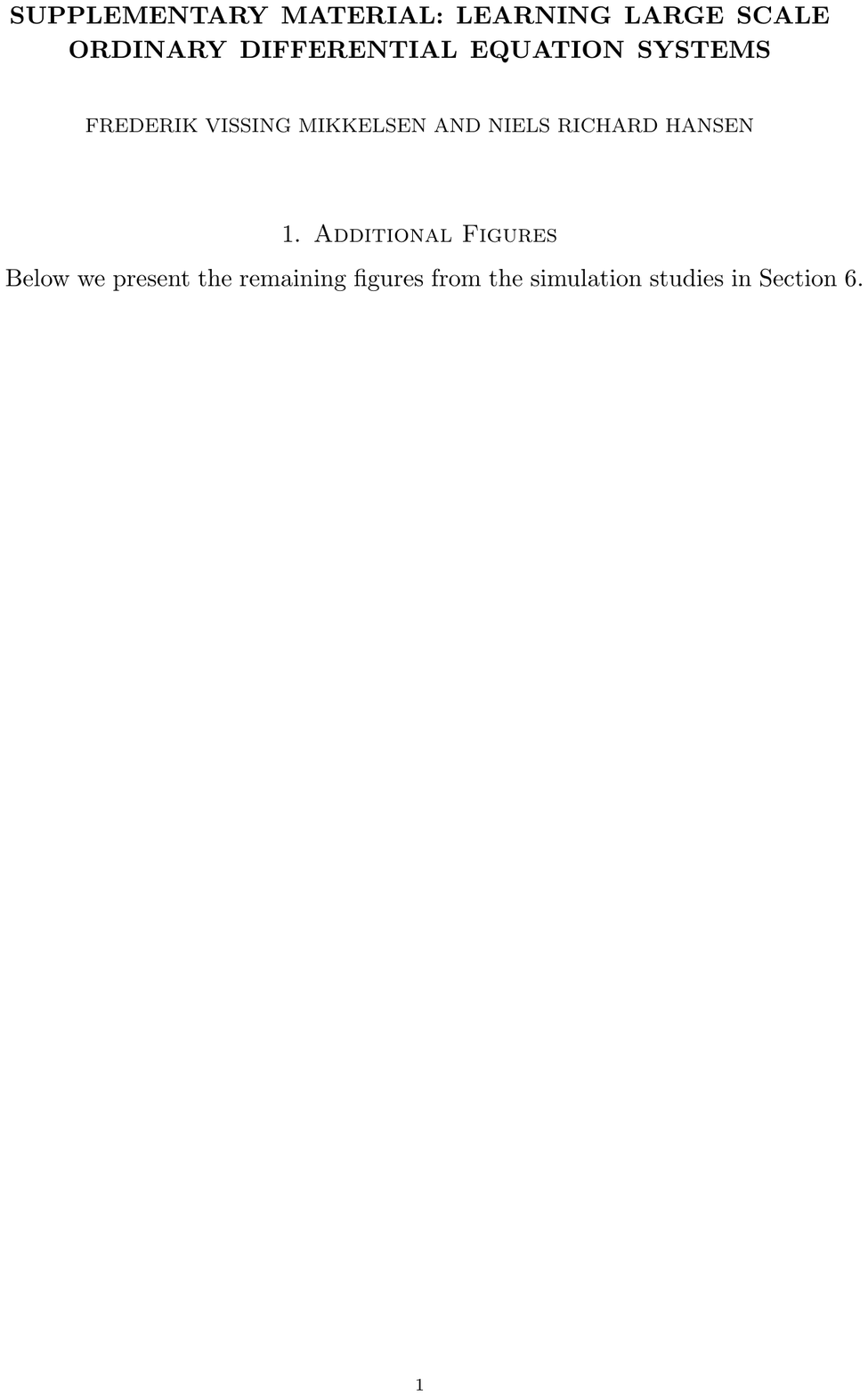}

\end{document}